\documentclass[11pt,fullpage]{article} 
\usepackage[centertags]{amsmath}
\usepackage{amsfonts}
\usepackage{amssymb}
\usepackage{amsthm}
\usepackage{eepic}
\usepackage{newlfont}
\usepackage{algorithm}
\usepackage{algorithmic}
\usepackage{enumitem}
\usepackage{tikz}
\usepackage[inactive]{srcltx}  
\theoremstyle{plain}
\newtheorem{thm}{Theorem}[section]
\newcommand{\BTHM}{\begin{thm}} \newcommand{\ETHM}{\end{thm}}
\newtheorem{cor}[thm]{Corollary}
\newcommand{\BCR}{\begin{cor}} \newcommand{\ECR}{\end{cor}}
\newtheorem{lem}[thm]{Lemma}
\newcommand{\BL}{\begin{lem}}   \newcommand{\EL}{\end{lem}}
\newtheorem{clm}[thm]{Claim}
\newcommand{\BCM}{\begin{clm}}   \newcommand{\ECM}{\end{clm}}
\newtheorem{prop}[thm]{Proposition}
\newcommand{\BP}{\begin{prop}}   \newcommand{\EP}{\end{prop}}
\newtheorem{assm}[thm]{Assumption}
\newcommand{\BASM}{\begin{assm}}   \newcommand{\EASM}{\end{assm}}

\theoremstyle{definition}
\newtheorem{defn}{Definition}[section]
\newcommand{\BD}{\begin{defn}}   \newcommand{\ED}{\end{defn}}
\newtheorem{con}[thm]{Conjecture}
\newcommand{\BCONJ}{\begin{con}}   \newcommand{\ECONJ}{\end{con}}

\theoremstyle{definition}
\newtheorem{problem}[thm]{Problem}
\newcommand{\BPR}{\begin{problem}}   \newcommand{\EPR}{\end{problem}}

\newenvironment{rem}{\noindent{\bf Remark:~~}}{}
\newcommand{\BREM}{\begin{rem}} \newcommand{\EREM}{\end{rem}}
\newenvironment{discussion}{\noindent{\bf Discussion:~~\\}}{}
\newcommand{\BDIS}{\begin{discussion}} \newcommand{\EDIS}{\end{discussion}}

\numberwithin{equation}{section}

\def\blackslug
     {\hbox{\hskip 1pt\vrule width 8pt height 8pt depth 1.5pt\hskip 1pt}}
\def\qed{\quad\blackslug\lower 8.5pt\null\par}

\newtheorem{exmp}[thm]{Example}
\newcommand{\BEX}{\begin{exmp}} \newcommand{\EEX}{\end{exmp}}
\newcommand{\BF}{\begin{fact}}   \newcommand{\EF}{\end{fact}}
\newcommand{\Bcr}{\begin{techcorr}}
\newcommand{\Ecr}{\end{techcorr}}
\newcommand{\BDS}{\begin{description}}
\newcommand{\EDS}{\end{description}}
\newcommand{\BE}{\begin{enumerate}}
\newcommand{\EE}{\end{enumerate}}
\newcommand{\BI}{\begin{itemize}}
\newcommand{\EI}{\end{itemize}}
\renewenvironment{proof}{\noindent{\bf Proof:~~}}{\qed}
\newcommand{\BPF}{\begin{proof}}
\newcommand{\EPF}{\end{proof}}

\newcommand{\BB}{\begin{enumerate}}
\newcommand{\EB}{\end{enumerate}}

\title{New results on large induced forests in graphs}

\author{Shimon Kogan \\ \\
  Department of Computer Science and Applied Mathematics \\
          Weizmann Institue, Rehovot 76100, Israel \\
          \ shimon.kogan@weizmann.ac.il
          }

\begin{document}
\maketitle
\begin{abstract}
For a graph $G$, let $a(G)$ denote the maximum size of a subset of vertices that induces a forest. We prove the following results.
\begin{enumerate}
  \item Let $G$ be a graph of order $n$, maximum degree $\Delta>0$ and maximum clique size $\omega$. Then
\[
a(G) \geq \frac{6n}{2\Delta + \omega +2}.
\]
This bound is sharp for cliques.
  \item Let $G=(V,E)$ be a triangle-free graph and let $d(v)$ denote the degree of $v \in V$. Then
  \[
  a(G) \geq \sum_{v \in V} \min\left(1, \frac{3}{d(v)+2} \right).
  \]
  As a corollary we have that a triangle-free graph $G$ of order $n$, with $m$ edges and average degree $d \geq 2$ satisfies
  \[
  a(G) \geq \frac{3n}{d+2}.
  \]
  This improves the lower bound $n - \frac{m}{4}$ of Alon-Mubayi-Thomas for graphs of average degree greater than $4$.
  Furthermore it improves the lower bound $\frac{20n - 5m - 5}{19}$ of Shi-Xu for (connected) graphs of average degree at least $\frac{9}{2}$.
\end{enumerate}
\end{abstract}




\section{Introduction}
For a (simple, undirected) graph $G = (V, E)$, we say that a set $S \subseteq V$ is an acyclic set if the
induced subgraph $G[S]$ is a forest. We let $a(G)$ denote the maximum size of an acyclic set in $G$.

In \cite{DBLP:journals/gc/AlonKS87} the following theorem was proven.
\BTHM\label{general_forest_bound1}
Let $G=(V,E)$ be a graph and let $d(v)$ denote the degree of $v \in V$. Then
  \[
  a(G) \geq \sum_{v \in V} \min\left(1, \frac{2}{d(v)+1} \right).
  \]
\ETHM
Furthermore the following corollary of Theorem \ref{general_forest_bound1} is shown in \cite{DBLP:journals/gc/AlonKS87}.
\BCR\label{general_forest_bound1_cor}
Let $G=(V,E)$ be a graph of order $n$ and average degree $d \geq 2$. Then
  \[
  a(G) \geq \frac{2n}{d+1}.
  \]
\ECR
In terms of maximum degree Corollary \ref{general_forest_bound1_cor} implies the following.
\BCR\label{general_forest_bound1_cor2}
Let $G=(V,E)$ be a graph of order $n$ and maximum degree $\Delta$. Then
  \[
  a(G) \geq \frac{2n}{\Delta+1}.
  \]
\ECR
A linear $k$-forest is a forest consisting of paths of length at most $k$ (that is the path contains $k$ edges). We let $a_k(G)$ denote the maximum size of an induced linear $k$-forest
in $G$.
We note that the following slight strengthening of Corollary \ref{general_forest_bound1_cor2} holds (See Appendix \ref{app1sect}).
\BTHM\label{general_forest_bound1_our_simple_thm}
Let $G=(V,E)$ be a graph of order $n$ and maximum degree $\Delta$. Then
  \[
  a_3(G) \geq \frac{2n}{\Delta+1}.
  \]
\ETHM
The following theorem was proven in \cite{DBLP:journals/jgt/AlonMT01}.
\BTHM\label{fundumental_bound_triangle_free_alon}
If $G$ is a triangle-free graph with $n$ vertices and $m$ edges, then $a(G) \geq n - \frac{m}{4}$.
\ETHM
Furthermore in \cite{DBLP:journals/jgt/ShiX17} the following is proven.
\BTHM\label{fundumental_bound_triangle_free_chinese}
If $G$ is a connected triangle-free graph with $n$ vertices and $m$ edges, then $a(G) \geq \frac{20n - 5m - 5}{19}$.
\ETHM
We note that Theorem \ref{fundumental_bound_triangle_free_chinese} has the following corollary (see Appendix \ref{last_appendix3}).
\BCR
Let $G=(V,E)$ be a triangle-free graph of order $n$ and average degree at most $4$. Then $a(G) \geq \frac{15n}{29}$.
\ECR
In this article we will prove the following theorem on triangle-free graphs.
\BTHM\label{our_bound_triangle_free_main}
Let $G=(V,E)$ be a triangle-free graph and let $d(v)$ denote the degree of $v \in V$. Then
  \[
  a(G) \geq \sum_{v \in V} \min\left(1, \frac{3}{d(v)+2} \right).
  \]
\ETHM
\BCR\label{our_bound_triangle_free_cor}
Let $G=(V,E)$ be a triangle-free graph of order $n$ and average degree $d \geq 2$. Then
  \[
  a(G) \geq \frac{3n}{d+2}.
  \]
\ECR
The bound in Corollary \ref{our_bound_triangle_free_cor} improves upon the bound in Theorem \ref{fundumental_bound_triangle_free_alon}
for graph of average degree greater than $4$. Furthermore the bound in Corollary \ref{our_bound_triangle_free_cor} improves upon the bound in Theorem
\ref{fundumental_bound_triangle_free_chinese} for graphs of average degree at least $\frac{9}{2}$.
Notice that the bound in Theorem \ref{fundumental_bound_triangle_free_chinese} holds only for connected triangle-free graphs while the bound in Corollary \ref{our_bound_triangle_free_cor}
holds for all triangle-free graphs (of average degree at least $2$). \\
We note that for large average degrees better bounds exist. 
In \cite{DBLP:journals/jct/AjtaiKS80} it is proved that every triangle-free graph on n vertices and average degree $d$ has
an independent set of size at least $\Omega\left(\frac{n \log d}{d}\right)$ (see also \cite{DBLP:journals/dm/Shearer83}).

The study of the size of a maximum acyclic set in graphs containing no clique of size $4$ was first addressed in \cite{DBLP:journals/jgt/AlonMT01}.
In particular the following theorem was proven in  \cite{DBLP:journals/jgt/AlonMT01}.
\BTHM\label{fundumental_bound_four_clique_free_alon}
If $G$ a graph with $n$ vertices and $m$ edges, such that $G$ contains no clique of size $4$ and $G$ has maximum degree $3$,
then $a(G) \geq n - \frac{m}{4} -\frac{1}{4}$.
\ETHM
We give general bounds on the size of a maximum acyclic set in terms of maximum degree and maximum clique size.
In particular we prove the following theorem.
\BTHM\label{beautiful_theorem1}
Let $G$ be a graph of order $n$, maximum degree $\Delta>0$ and maximum clique size $\omega$. Then
\[
a(G) \geq \frac{6n}{2\Delta + \omega +2}.
\]
\ETHM
Notice that this bound is sharp for cliques. Furthermore the forest obtained in Theorem \ref{beautiful_theorem1} is in fact linear in the case of $\omega\geq 4$.

The theorem above is an analogue of the following theorem on independent sets which is proven in \cite{Fajtlowicz1}.
\BTHM\label{faj1}
Let $G$ be a graph of order $n$, maximum degree $\Delta$ and maximum clique size $\omega$. Let $\alpha(G)$ denote the size of the maximum independent
set of $G$. Then
\[
\alpha(G) \geq \frac{2n}{\Delta + \omega + 1}.
\]
\ETHM
We mention the following result which was proven in \cite{kostochka1982modification} and \cite{DBLP:journals/jgt/Rabern13}.
\BTHM
The vertex set of any triangle-free graph $G$ of maximum degree $\Delta$ can be partitioned into $\lceil \frac{\Delta+2}{3} \rceil$ sets, each of which induces
a disjoint union of paths in $G$.
\ETHM
This result is in a sense complementary to Corollary \ref{our_bound_triangle_free_cor}. \\
Another result related to this paper is Theorem $6.1$ of \cite{DBLP:journals/jgt/CranstonR15} (first proven in \cite{borodin1976decomposition}).
We state only a special case of this theorem related to forests.
\BTHM
Let $G$ be a graph of maximum degree $\Delta \geq 4$ containing no cliques of size $\Delta+1$. Then the vertex set of graph $G$ can be partitioned into
$\lceil \frac{\Delta}{2} \rceil$ sets, each of which induces a disjoint union of paths in $G$.
\ETHM

\section{Triangle-free graphs}
In this section we will prove Theorem \ref{our_bound_triangle_free_main} and Corollary \ref{our_bound_triangle_free_cor}. \\
Define the potential function $f(d) = \min\left(1, \frac{3}{d+2} \right)$.
We shall require two technical lemmas.
\BL\label{triangle_free_technical_lem1}
Let $\Delta \geq 5$ be an integer. Let $2 \leq d \leq \Delta$ and $0 \leq q \leq d$  be integers. Then
\[
f(d-q) - f(d) \geq  q ( f(\Delta-1) - f(\Delta) ).
\]
\EL
\BPF
If $q<d$ the claim follows from the inequality $f(d-1) - f(d) \leq f(d-2) - f(d-1)$ which holds for all $d \geq 3$ and the fact that
\[
f(d-q) - f(d) = ((f(d-1)-f(d)) + (f(d-2)-f(d-1)) + \ldots + (f(d-q)-f(d-q+1))
\]
If $q=d$ then $f(d - q) = 1$ and we need to prove that
\begin{equation}\label{accounting_lemma_1_triangle_free}
1 - f(d) \geq d( f(\Delta-1) - f(\Delta) ).
\end{equation}
Inequality \ref{accounting_lemma_1_triangle_free} holds if and only if
\begin{equation}\label{accounting_lemma_2_triangle_free}
1 - \frac{3}{d+2} \geq d \left( \frac{3}{\Delta+1} - \frac{3}{\Delta+2} \right).
\end{equation}
As $2 \leq d \leq \Delta$ and $f(d-1) - f(d) \leq f(d-2) - f(d-1)$ for $d \geq 3$, Inequality \ref{accounting_lemma_2_triangle_free} holds if
\begin{equation}\label{accounting_lemma_3_triangle_free}
1 - \frac{3}{d+2} \geq d\left( \frac{3}{d+1} - \frac{3}{d+2} \right).
\end{equation}
And Inequality \ref{accounting_lemma_3_triangle_free} holds if and only if
\begin{equation}\label{accounting_lemma_4_triangle_fre}
\frac{3d}{d+1} \leq d-1
\end{equation}
This inequality holds for $d \geq 4$.
Hence we may assume that $2 \leq d \leq 3$. As $\Delta \geq 5$ we can verify that
\begin{equation}\label{accounting_lemma_4_triangle_free2}
1 - \frac{3}{d+2} \geq d\left( \frac{3}{6} - \frac{3}{7} \right) \geq d\left( \frac{3}{\Delta+1} - \frac{3}{\Delta+2}\right).
\end{equation}
And thus Inequality \ref{accounting_lemma_2_triangle_free} follows and we are done.
\EPF
The following lemma is almost identical to Lemma \ref{triangle_free_technical_lem1}, we give a proof for completeness.
\BL\label{triangle_free_technical_lem2}
Let $\Delta \geq 5$ be an integer. Let $2 \leq d < \Delta$ and $0 \leq q \leq d$  be integers. Then
\[
f(d-q) - f(d) \geq  q ( f(\Delta-2) - f(\Delta-1) ).
\]
\EL
\BPF
If $q<d$ the claim follows immediately from the inequality $f(d-1) - f(d) \leq f(d-2) - f(d-1)$ which holds for all $d \geq 3$. \\
If $q=d$ then $f(d - q) = 1$ and we need to prove that
\begin{equation}\label{accounting_lemma_1_triangle_free2}
1 - f(d) \geq d( f(\Delta-2) - f(\Delta-1) ).
\end{equation}
Inequality \ref{accounting_lemma_1_triangle_free2} holds if and only if
\begin{equation}\label{accounting_lemma_2_triangle_free2}
1 - \frac{3}{d+2} \geq d\left( \frac{3}{\Delta} - \frac{3}{\Delta+1} \right).
\end{equation}
As  $2 \leq d \leq \Delta-1$ and $f(d-1) - f(d) \leq f(d-2) - f(d-1)$ for $d \geq 3$, Inequality \ref{accounting_lemma_2_triangle_free2} holds if
\begin{equation}\label{accounting_lemma_3_triangle_free2}
1 - \frac{3}{d+2} \geq d\left( \frac{3}{d+1} - \frac{3}{d+2} \right).
\end{equation}
And we have shown in Lemma \ref{triangle_free_technical_lem1} (Inequality \ref{accounting_lemma_3_triangle_free})
that Inequality \ref{accounting_lemma_3_triangle_free2} holds
for $d \geq 4$.
Hence we may assume that $2 \leq d \leq 3$. As $\Delta \geq 5$ we can verify that
\begin{equation}\label{accounting_lemma_5_triangle_free2}
1 - \frac{3}{d+2} \geq d\left( \frac{3}{5} - \frac{3}{6} \right) \geq d\left( \frac{3}{\Delta} - \frac{3}{\Delta+1}\right).
\end{equation}
And thus Inequality \ref{accounting_lemma_2_triangle_free2} follows and we are done.
\EPF
\text{} \\
\noindent\textbf{Proof of Theorem \ref{our_bound_triangle_free_main}}: \\
Let $G=(V,E)$ be a triangle-free graph and let $d(v)$ denote the degree of $v \in V$. We shall prove that
  \[
  a(G) \geq \sum_{v \in V} \min\left(1, \frac{3}{d(v)+2} \right).
  \]
Let $n$ be the number of vertices in graph $G$ and $m$ be the number of edges in graph $G$.
Let $\Delta$ be the maximum degree of graph $G$ and $\delta$ the minimum degree of graph $G$.
We shall prove the theorem by induction on $n$. Clearly it holds for $n=1$.
Suppose that the vertices of $G$ are $v_1,\ldots,v_n$ and that the degree sequence of graph $G$ is $d_1,\ldots,d_n$.
Let $S= \sum_{i=1}^{n} f(d_i)$.
We need to prove that $a(G) \geq S$. \\
First Assume that there is a vertex in $G$ of degree at most $1$. \\
Assume without loss of generality this vertex is $v_1$, that is $d_1 \leq 1$.
Define $H$ to be the graph formed from $G$ by removing vertex $v_1$  from $G$. Let $d_1', \ldots d_{n-1}'$ be the degree
sequence of graph $H$. Let $T = \sum_{i=1}^{n-1} f(d_i')$. Notice that $T \geq S-1$. By the induction hypothesis $a(H) \geq T$.
Furthermore we can add $v_1$ to a maximum forest in $H$ and the resulting set will be a forest in $G$.
Hence $a(G) \geq a(H) + 1 \geq T+1 \geq S$ and we are done. \\
Henceforth we shall assume that all the vertices of $G$ are of degree at least $2$, that is $\delta \geq 2$. \\
Assume that $\Delta \leq 4$. By Theorem \ref{fundumental_bound_triangle_free_alon} we have
\[
a(G) \geq n - \frac{m}{4} = \sum_{i=1}^{n} \left(1-\frac{d_i}{8}\right).
\]
As $2 \leq d_i \leq 4$ for all $i$, we have that $1-\frac{d_i}{8} \geq \frac{3}{d_i+2}$ for all $i$. And thus $a(G) \geq S$ and we are done. \\
Henceforth we assume that $\Delta \geq 5$ and $\delta \geq 2$. \\
We choose a vertex $v$ in $G$ such that the following conditions are satisfied.
\begin{enumerate}[label=(\arabic*), leftmargin=1.5cm]
  \item $d(v) = \Delta$.
  \item Subject to $(1)$, the number of neighbors of $v$ of degree $\Delta$ is maximized.
\end{enumerate}
Assume w.l.o.g that vertex $v$ chosen in the process above is $v_1$.
Let $T$ be the set of indices of the neighbors of $v_1$ of degree $\Delta$. Let $t$ be the number of neighbors of $v_1$ of degree $\Delta$, that is  $t=|T|$.
Define $S_1$ to be the set of indices of the neighbors of $v_1$ . Define $S_2$ to be the set of indices of the vertices which are at distance $2$ from $v_1$.
We consider $3$ cases. \\
\textbf{Case 1: $t=0$.} \\
Define $H$ to be the graph formed from $G$ by removing vertex $v_1$. Let $d'_1, \ldots d'_{n-1}$ be the degree sequence of graph $H$. Let
$Q = \sum_{i=1}^{n-1} f(d'_i)$. Now notice that
\begin{align}\label{triange_free_main_thm_case1}
Q &= S - f(d_1) + \sum_{i \in S_1} (f(d_i - 1) -f(d_i)) \notag \\
&\geq S - \frac{3}{\Delta+2} + \Delta \left( \frac{3}{\Delta} - \frac{3}{\Delta+1}\right) \geq S. \notag
\end{align}
And we are done by applying the induction hypothesis to graph $H$. \\
\textbf{Case 2: $t=\Delta$.} \\
Define $T' \subseteq T$ to be the set of indices of arbitrary $\Delta-1$ neighbors of $v_1$.
For every $i$ denote by $n_i$ the number of neighbors of vertex $v_i$ with indices in $T'$. \\
Define $H$ to be the graph formed from $G$ by removing the vertices with indices in $T'$. Let $d'_1, \ldots d'_{n-\Delta+1}$ be the degree sequence of graph $H$.
Let $Q = \sum_{i=1}^{n-\Delta+1} f(d'_i)$. Now notice that
\begin{equation}\label{triangle_free1_full_sum1}
Q = S + (1-f(\Delta))  - (\Delta-1)f(\Delta) + \sum_{i \in S_2} (f(d_i-n_i) - f(d_i)).
\end{equation}
Notice that in Equation \ref{triangle_free1_full_sum1}, $(1-f(\Delta))$ is the change of potential for vertex $v_1$. And $-(\Delta-1)f(\Delta)$ is the potential change
from deleting the vertices with indices in $T'$. Finally $\sum_{i \in S_2} (f(d_i-n_i) - f(d_i))$ is the potential change to vertices with indices in $S_2$.\\
By Lemma \ref{triangle_free_technical_lem1} we have
\begin{equation}\label{triange_free_lemlemlem1}
f(d_i-n_i) - f(d_i) \geq n_i( f(\Delta-1) - f(\Delta) ).
\end{equation}
And thus from Equations \ref{triangle_free1_full_sum1} and \ref{triange_free_lemlemlem1} we have
\begin{equation}\label{triangle_free1_full_sum2}
Q \geq S + (1-f(\Delta)) - (\Delta-1)f(\Delta) + \sum_{i \in S_2} n_i(f(\Delta-1) - f(\Delta)).
\end{equation}
Now notice that since each vertex with an index in $T'$ has $\Delta-1$ neighbors with indices in $S_2$ we have
\begin{equation}\label{triange_free_degree_consideration1}
\sum_{i \in S_2} n_i = (\Delta-1)^2.
\end{equation}
Hence from Equations \ref{triangle_free1_full_sum2} and \ref{triange_free_degree_consideration1} we have
\begin{align}\label{triange_free_main_thm_case2}
  Q & \geq  S + (1-f(\Delta) - (\Delta-1)f(\Delta) + (\Delta-1)^2(f(\Delta-1) - f(\Delta)) \notag \\
   & = S +\left(1-\frac{3}{\Delta+2}\right)- (\Delta-1) \frac{3}{\Delta+2} + (\Delta-1)^2 \left( \frac{3}{\Delta+1} - \frac{3}{\Delta+2} \right)\notag  \\
   & = S + \frac{ (\Delta -1)(\Delta-5) } { (\Delta +1)(\Delta+2)  } \geq S &&\text{( as $\Delta \geq 5$ ) }\notag
\end{align}
And we are done by applying the induction hypothesis to graph $H$. \\
\textbf{Case 3: $0< t <\Delta$.} \\
For every $i$ denote by $n_i$ the number of neighbors of vertex $v_i$ with indices in $T$. \\
Define $H$ to be the graph formed from $G$ by removing the vertices with indices in $T$. Let $d'_1, \ldots d'_{n-t}$ be the degree sequence of graph $H$.
Let $Q = \sum_{i=1}^{n-t} f(d'_i)$. Now notice that
\begin{equation}\label{triangle_free4_full_sum1}
Q = S + (f(\Delta-t)-f(\Delta))  - t f(\Delta) + \sum_{i \in S_2} (f(d_i-n_i) - f(d_i)).
\end{equation}
Notice that in Equation \ref{triangle_free4_full_sum1}, $(f(\Delta-t)-f(\Delta))$ is the change of potential for vertex $v_1$. And $-t f(\Delta)$ is the potential change
from deleting the vertices with indices in $T$. \\
Finally $\sum_{i \in S_2} (f(d_i-n_i) - f(d_i))$ is the potential change to vertices with indices in $S_2$.\\
Let $A \subseteq S_2$ be the set of indices in $S_2$ of vertices of degree $\Delta$ in $G$.
Let $B \subseteq S_2$ bet set of indices in $S_2$ of vertices of degree at most $\Delta-1$ in $G$. Notice that $S_2 = A \cup B$. \\
Hence we may rewrite Equation \ref{triangle_free4_full_sum1} as
\begin{equation}\label{triangle_free7_full_sum2}
Q = S + (f(\Delta-t)-f(\Delta))  - t f(\Delta) + \sum_{i \in A} (f(d_i-n_i) - f(d_i)) +  \sum_{i \in B} (f(d_i-n_i) - f(d_i)).
\end{equation}
By Lemma \ref{triangle_free_technical_lem1} we have that for all $i \in A$ the following holds.
\begin{equation}\label{triange_free7_lemlemlem1}
f(d_i-n_i) - f(d_i) \geq n_i( f(\Delta-1) - f(\Delta) ).
\end{equation}
Furthermore by Lemma \ref{triangle_free_technical_lem2} we have that for all $i \in B$ the following holds.
\begin{equation}\label{triange_free7_lemlemlem2}
f(d_i-n_i) - f(d_i) \geq n_i( f(\Delta-2) - f(\Delta-1) ).
\end{equation}
Applying Inequalities \ref{triange_free7_lemlemlem1} and \ref{triange_free7_lemlemlem2} to Equation \ref{triangle_free7_full_sum2} we get
\begin{equation}\label{triangle_free7_full_sum3}
Q \geq  S + (f(\Delta-t)-f(\Delta))  - t f(\Delta) + \sum_{i \in A} n_i (f(\Delta-1) - f(\Delta)) +  \sum_{i \in B} n_i (f(\Delta-2) - f(\Delta-1)).
\end{equation}
As each vertex with an index in $T$ has at most $t-1$ neighbors of degree $\Delta$ besides $v_1$
(this follows from condition $(2)$) we have the following inequality.
\begin{equation}\label{triangle_free_crux1}
\sum_{i \in A} n_i \leq t(t-1).
\end{equation}
Now as
\begin{equation}\label{triangle_free_crux2}
\sum_{i \in S_2} n_i = t(\Delta-1).
\end{equation}
We conclude that
\begin{equation}\label{triangle_free_crux3}
\sum_{i \in B} n_i \geq t(\Delta-1) - t(t-1) = t(\Delta-t).
\end{equation}
Applying Inequalities \ref{triangle_free_crux1} and \ref{triangle_free_crux3} to Equation \ref{triangle_free7_full_sum3} we get
\begin{align}
Q - S & \geq  f(\Delta-t)-f(\Delta)  - t f(\Delta) + t(t-1) (f(\Delta-1) - f(\Delta)) +  t(\Delta-t) (f(\Delta-2) - f(\Delta-1)) \notag \\
& = \frac{3}{\Delta-t+2} - \frac{3}{\Delta+2} - t \frac{3}{\Delta+2}+t(t-1)\left( \frac{3}{\Delta+1} - \frac{3}{\Delta+2} \right)
+ t(\Delta-t) \left( \frac{3}{\Delta} - \frac{3}{\Delta+1} \right)  \notag \\
& = \frac{3t( \Delta^2 - 2\Delta t + \Delta + 2t^2 - 4t  )}{\Delta(\Delta+1)(\Delta+2)(\Delta-t + 2)} \notag \\
&= \frac{3t[ (\Delta - t)^2 + (t-2)^2 + \Delta-4  ]}{\Delta(\Delta+1)(\Delta+2)(\Delta-t + 2)} \geq 0 \text{ \hspace{6cm} ( as $\Delta \geq 5$ ) } \notag
\end{align}
And we are done by applying the induction hypothesis to graph $H$.
\qed \text{} \\
\noindent\textbf{Proof of Corollary \ref{our_bound_triangle_free_cor}}: \\
Let $G=(V,E)$ be a triangle-free graph of order $n$, with $m$ edges and average degree $d \geq 2$. We shall prove that
  \[
  a(G) \geq \frac{3n}{d+2}.
  \]
This proof is similar to the proof of Corollary $1.4$ in \cite{DBLP:journals/gc/AlonKS87}. \\
Recall that by Theorem \ref{fundumental_bound_triangle_free_alon} we have
$a(G) \geq n - m/4 = n(1-d/8)$. Now notice that for $2 \leq d \leq 4$ we have that $n\left(1- \frac{d}{8}\right) \geq \frac{3n}{d+2}$.
Hence $a(G) \geq \frac{3n}{d+2}$ for $2 \leq d \leq 4$. Henceforth we assume that $d>4$.
By Theorem \ref{our_bound_triangle_free_main} we have that $a(G) \geq w$, where $w$ is the minimum possible value of the expression
\begin{equation}\label{alon_trick1}
\sum_{i=1}^n \min\left(1, \frac{3}{q_i+2} \right)
\end{equation}
subject to the constraints
\begin{equation}\label{alon_trick2}
\sum_{i=1}^n q_i = 2m  \text{\quad and \quad} q_i \geq 0 \text{ are integers }
\end{equation}
If there is an $i$ such that $q_i = 0$ then as $d>4$ we have some $j$ for which $q_j = r \geq 5$.
Setting $q_i = 2$ and $q_j = r - 2$, we get a new sequence which
decreases the sum of \ref{alon_trick1} (as $\frac{1}{4} > \frac{3}{r} - \frac{3}{r+2}$ for $r\geq 5$) thus contradicting the minimality of w.
Hence we may assume that for all $1 \leq i \leq n$, $q_i \geq 1$. We conclude that
\[
a(G) \geq \sum_{i=1}^n \frac{3}{q_i+2} \geq \frac{3n}{d+2}
\]
where the last inequality follows from Jensen's inequality. \qed

\section{Graphs without a clique of size $4$}\label{four_clique_section}
In this section we shall prove the following theorem.
\BTHM\label{four_clique_main_thm}
Let $G=(V,E)$ be a graph of order $n$ and maximum degree $\Delta>0$, containing no cliques of size $4$. Then
  \[
  a(G) \geq \frac{6n}{2\Delta+5}.
  \]
\ETHM
We shall start by proving the following lemma.
\BL\label{four_clique_main_lemma}
Let $\Delta>0$ and let $G=(V,E)$ be a $\Delta$-regular graph of order $n$, containing no cliques of size $4$. Then
  \[
  a(G) \geq \frac{6n}{2\Delta+5}.
  \]
\EL
\BPF
Given a set $S$ of vertices of $G$, Let $|S|$ denote the number of vertices in $S$, $G[S]$ denote the subgraph of $G$ induced by the vertices of $S$, and
$e(S)$ denote the number of edges in $G[S]$. Given an induced subgraph $T$ of $G[S]$ we denote by $\Delta(T)$ the maximum degree of $T$.
We denote by $D(T)$ the diameter of $T$ (that is the greatest distance between any pair of vertices in $T$).
Finally we denote by $P(T)$ the number of paths in $T$ of length $D(T)$. \\
Choose an induced forest $S$ in graph $G$ such that the following conditions are satisfied.
\begin{enumerate}[label=(\arabic*), leftmargin=1.5cm]
  \item $|S|$ is maximized.
  \item Subject to $(1)$, $e(S)$ is maximized.
  \item Subject to $(2)$, the number of vertices of degree $1$ in $G[S]$ is maximized.
  \item  Subject to $(3)$, we maximize the following sum.
  \[
  \sum_{T \text{ is a tree in } G[S]} \Delta(T).
  \]
  \item Subject to $(4)$, we minimize the following sum.
  \[
  \sum_{T \text{ is a tree in } G[S]} P(T).
  \]
\end{enumerate}
Let $\beta_i(S)$ denote the number of vertices in $V \backslash S$ with exactly $i$ adjacent vertices in $S$.
Notice that
\begin{equation}\label{basic_step_1_main4}
  n - |S| = \sum_{i=0}^{|S|} \beta_i(S).
\end{equation}
Given a vertex $s \in S$, let $d^o(s)$ denote the number of vertices in $V \backslash S$ that are adjacent to $s$.
Notice that
\begin{equation}\label{basic_step_2_main4}
 \sum_{s \in S} d^o(s) = \sum_{i=0}^{|S|} i \beta_i(S).
\end{equation}
Multiplying Equation \ref{basic_step_1_main4} by $3$ and subtracting Equation \ref{basic_step_2_main4} we obtain the following.
\begin{equation}\label{basic_step_3_main4}
3n - 3|S| - \sum_{s \in S} d^o(s) =  \sum_{i=0}^{|S|} (3 - i) \beta_i(S).
\end{equation}
By condition $(1)$ we have $\beta_1(S) = \beta_0(S) = 0$.  Thus from Equation \ref{basic_step_3_main4} it follows that
\begin{equation}\label{basic_step_4_main4}
3n - 3|S| - \sum_{s \in S} d^o(s) = \beta_2(S) - \sum_{i=4}^{|S|} (i - 3) \beta_i(S).
\end{equation}
We conclude that
\begin{equation}\label{basic_step_5_main4}
3n - 3|S| - \sum_{s \in S} d^o(s) \leq \beta_2(S) - \sum_{i=4}^{|S|} \beta_i(S)
\end{equation}
Notice that
\begin{equation}\label{basic_step_6_main4}
\sum_{s \in S} d^o(s) \leq \Delta |S| - \sum_{T \text{ is a tree in } G[S]} 2e(T).
\end{equation}
Where $e(T)$ is the number of edges in tree $T$.
Hence by Equations \ref{basic_step_5_main4} and \ref{basic_step_6_main4} we have
\begin{equation}\label{basic_step_7_main4}
3n - 3|S| - \Delta |S| \leq  \beta_2(S) -
\sum_{T \text{ is a tree in } G[S]} 2e(T) - \sum_{i=4}^{|S|} \beta_i(S).
\end{equation}
Hence if the following Inequality is satisfied (for $\Delta>0$)
\begin{equation}\label{basic_step_8_main4}
\sum_{T \text{ is a tree in } G[S]} 2e(T) - \beta_2(S) + \sum_{i=4}^{|S|} \beta_i(S) \geq \frac{|S|}{2}.
\end{equation}
Then we are done as from Inequalities \ref{basic_step_8_main4} and \ref{basic_step_7_main4} we have
\[
3n - 3|S| - \Delta |S| \leq -\frac{|S|}{2}.
\]
And thus
\[
|S| \geq \frac{6n}{2\Delta+5}.
\]
The rest of the section is devoted to the proof of Inequality \ref{basic_step_8_main4}. \\
Let Q be the set of vertices in $V \backslash S$ with at least $4$ adjacent vertices in $S$.
As $|Q| = \sum_{i=4}^{|S|} \beta_i(S)$ we need to prove that
\begin{equation}\label{basic_step_9_main4}
\sum_{T \text{ is a tree in } G[S]} 2e(T) - \beta_2(S) + |Q| \geq \frac{|S|}{2}.
\end{equation}
Let $S_0$ denote the set of vertices of degree $0$ in $G[S]$. Let $B_i$ denote the set of vertices of $V \backslash S$ with exactly $i$ adjacent vertices in $S$.
Notice that by definition $|B_i| = \beta_i$.
Given a vertex $s$ in $S$ we denote by $d_S(s)$ the degree of vertex $s$ in $G[S]$.
We shall need the following observations.
\\ \text{} \\
\noindent\textbf{Observation 1:} No vertex in $B_2$ is adjacent to a vertex in $S_0$. \\
\textbf{Proof:} If such vertex $v \in B_2$ exists we can add it to $S$ and get a contradiction to condition $(1)$. \qed \text{} \\
\noindent\textbf{Observation 2:} Any vertex $s \in S$ is adjacent to at most two vertices in $B_2$. \\
\textbf{Proof:} Assume by contradiction that vertex $s$ is adjacent to vertices $v_1,v_2,v_3$ in $B_2$. As graph $G$ has no cliques of size $4$ we may assume without loss of generality
that vertices $v_1$ and $v_2$ are not adjacent. We remove vertex $s$ from $S$ and add vertices $v_1$ and $v_2$ to $S$, thus getting a contradiction to condition $(1)$. \qed \text{} \\
\noindent\textbf{Observation 3:} Any tree $T$ in $G[S]$ has at most $|T|$ adjacent vertices in $B_2$. \\
\textbf{Proof:} Each vertex in $B_2$ is either adjacent to two vertices in $T$ or not adjacent to any vertex in $T$, for otherwise we get a contradiction to condition $(1)$.
Hence Observation $3$ follows from Observation $2$ by double counting
(as there are at most $2|T|$ edges between $T$ and $B_2$ and each vertex in $B_2$ that is adjacent to a vertex in $T$ must be adjacent to exactly two vertices in $T$)
. \qed \text{} \\
\noindent\textbf{Observation 4:} Any tree $T$ in $G[S]$ for which $|T| \leq 7$ has at most $|T|-1$ adjacent vertices in $B_2$. \\
\textbf{Proof:} Given a tree $T$ on at most $7$ vertices, we know by Observation $3$ that $T$ has at most $|T|$ adjacent vertices in $B_2$.
We shall show that in fact $T$ has at most $|T|-1$ adjacent vertices in $B_2$.
We shall do a case analysis on all non-isomorphic trees of at most $7$ vertices. Due to the length of the case analysis we shall prove this claim in Appendix \ref{app2sect}.
 \qed \text{} \\
\noindent\textbf{Observation 5:} If a vertex $v \in B_3$ is adjacent to a vertex in $S_0$ then the two other neighbors of $v$ in $S$, which we denote by $s_1$ and $s_2$, satisfy
$d_S(s_1) \geq 2$ and $d_S(s_2) \geq 2$. \\
\textbf{Proof:} We prove by contradiction. Assume without loss of generality that $d_S(s_1) \leq 1$. We remove vertex $s_1$ from $S$ and add vertex $v$, thus getting a contradiction to
condition $(2)$. \qed \text{} \\
\noindent\textbf{Observation 6:} A vertex $v \in B_3$ can be adjacent to at most one vertex in $S_0$. \\
\textbf{Proof:}
Follows from Observation $5$. \qed \text{} \\
Let $T_1,\ldots,T_t$ be the trees in $G[S]$ such that each such tree has at least $2$ vertices and at most $7$ vertices.
Let $T_{t+1},T_{t+2},\ldots,T_k$ be the trees in $G[S]$ of at least $8$ vertices.
By Observation $3$ and Observation $4$ (and the fact that tree $T_i$ had $|T_i|-1$ edges) we have
\begin{equation}\label{basic_step_10_main4}
\begin{split}
\sum_{T \text{ is a tree in } G[S]} 2e(T) - \beta_2(S) & \geq \sum_{i=1}^{t} (|T_i|-1) + \sum_{i=t+1}^{k} (|T_i| -2) \\
 & = |S| - |S_0| - t - 2(k-t).
\end{split}
\end{equation}
As $\sum_{i=1}^{t} |T_i| \geq 2t$ we have $\sum_{i=t+1}^{k} |T_i| \leq |S|-|S_0|-2t$, and thus
\begin{equation}\label{basic_step_11_main4}
k-t \leq \frac{|S|-|S_0|-2t}{8}.
\end{equation}
Combining Equations \ref{basic_step_10_main4} and \ref{basic_step_11_main4} we get
\begin{equation}\label{basic_step_12_main4}
\sum_{T \text{ is a tree in } G[S]} 2e(T) - \beta_2(S) \geq \frac{3|S|-3|S_0|-2t}{4}.
\end{equation}
As graph $G$ is $\Delta$-regular we have by Observation $1$ and Observation $6$ that the set $S_0$ has at least
$\Delta(|S_0|-|Q|)$ adjacent vertices in $B_3$. Hence by Observation $5$ there are at least $2(|S_0|-|Q|)$ vertices of degree at least $2$ in $G[S]$.
Hence we have in $G[S]$
\begin{itemize}
  \item $|S_0|$ vertices of degree $0$.
  \item at least $2t$ vertices of degree $1$ (leaves of trees).
  \item at least $2(|S_0|-|Q|)$ vertices of degree at least $2$.
\end{itemize}
We conclude that
\begin{equation}\label{basic_step_13_main4}
|S_0| + 2t + 2(|S_0|-|Q|) \leq |S|.
\end{equation}
And hence
\begin{equation}\label{basic_step_14_main4}
3|S_0| + 2t \leq |S| + 2|Q|.
\end{equation}
Combining Equation \ref{basic_step_12_main4} and Equation \ref{basic_step_14_main4} we get
\begin{equation}\label{basic_step_15_main4}
\sum_{T \text{ is a tree in } G[S]} 2e(T) - \beta_2(S) +|Q| \geq \frac{|S| + |Q|}{2} \geq \frac{|S|}{2}.
\end{equation}
And thus we are done.
\EPF
We shall prove now that Theorem \ref{four_clique_main_thm} follows from Lemma \ref{four_clique_main_lemma}. \\
Observe that, if $G$ is a graph with maximum degree $\Delta>0$, then we can create a $\Delta$-regular
graph by taking copies $H_1,H_2, \dots , H_r$ of $G$ and joining some pairs of vertices from different copies so
as to make the resulting graph $G'$ a $\Delta$-regular graph. This can be done without creating cliques of size $4$ if sufficiently many copies of $G$ are used.
Applying Lemma \ref{four_clique_main_lemma} to graph $G'$
we get by the pigeonhole principle that for some $1 \leq i \leq r$ we have $a(H_i) \geq \frac{6n}{2\Delta+5}$ and thus we are done.

\section{Graphs without a clique of size $q \geq 5$}
Recall that a linear $k$-forest is a forest consisting of paths of length at most $k$. We let $a_k(G)$ denote the maximum size of an induced linear $k$-forest
in $G$.
\BTHM\label{five_theorem_main}
Let $G=(V,E)$ be a graph of order $n$ and maximum degree $\Delta$, containing no cliques of size $q \geq 5$. Then
  \[
  a_4(G) \geq \frac{6n}{2\Delta+q+1}.
  \]
\ETHM
\BPF
Given a set $S$ of vertices of $G$, Let $|S|$ denote the number of vertices in $S$, $G[S]$ denote the subgraph of $G$ induced by the vertices of $S$, and
$e(S)$ denote the number of edges in $G[S]$.
Choose an induced linear $4$-forest $S$ in graph $G$ such that the following conditions are satisfied.
\begin{enumerate}[label=(\arabic*), leftmargin=1.5cm]
  \item $3|S| - e(S)$ is maximized.
  \item Subject to $(1)$, the number of vertices of degree $0$ in $G[S]$ is maximized.
\end{enumerate}
Let $\beta_i(S)$ denote the number of vertices in $V \backslash S$ with exactly $i$ adjacent vertices in $S$.
Notice that
\begin{equation}\label{basic_step_1_main5}
  n - |S| = \sum_{i=0}^{|S|} \beta_i(S).
\end{equation}
Given a vertex $s \in S$, let $d^o(s)$ denote the number of vertices in $V \backslash S$ that are adjacent to $s$.
Notice that
\begin{equation}\label{basic_step_2_main5}
 \sum_{s \in S} d^o(s) = \sum_{i=0}^{|S|} i \beta_i(S).
\end{equation}
Multiplying Equation \ref{basic_step_1_main5} by $3$ and subtracting Equation \ref{basic_step_2_main5} we obtain the following.
\begin{equation}\label{basic_step_3_main5}
3n - 3|S| - \sum_{s \in S} d^o(s) =  \sum_{i=0}^{|S|} (3 - i) \beta_i(S).
\end{equation}
By condition $(1)$ we have $\beta_0(S) = 0$.  Now we shall prove that $\beta_1(S) = 0$. \\
Assume by contradiction that there is a vertex $v$ in $V \backslash S$ with exactly one adjacent vertex $s$ in $S$.
If $s$ is a vertex of degree at least $1$ in $G[S]$ then we remove vertex $s$ from $S$ and add vertex $v$ to $S$, thus getting a contradiction to condition $(1)$.
Hence $s$ is a vertex of degree $0$ in $G[S]$. Thus we can add $v$ to $S$ and get a contradiction to condition $(1)$. \\
We have shown that $\beta_0(S) = 0$ and $\beta_1(S) = 0$. Thus from Equation \ref{basic_step_3_main5} it follows that
\begin{equation}\label{basic_step_4_main5}
3n - 3|S| - \sum_{s \in S} d^o(s) \leq \beta_2(S)
\end{equation}
Given a path $P$ we denote by $l(P)$ the length of the path (that is the number of edges in the path).
Notice that
\begin{equation}\label{basic_step_5_main5}
\sum_{s \in S} d^o(s) \leq \Delta |S| - \sum_{P \text{ is a path in } G[S]} 2l(P)
\end{equation}
Hence it follows from Equations \ref{basic_step_5_main5} and \ref{basic_step_4_main5} that
\begin{equation}\label{basic_step_6_main5}
3n - 3|S| - \Delta|S| \leq \beta_2(S) - \sum_{P \text{ is a path in } G[S]} 2l(P)
\end{equation}
Subtracting $\frac{q-5}{2}|S|$ from both sides we get
\begin{equation}\label{basic_step_7_main5}
3n - 3|S| - \Delta|S| - \frac{q-5}{2}|S| \leq \beta_2(S) - \frac{q-5}{2}|S| - \sum_{P \text{ is a path in } G[S]} 2l(P)
\end{equation}
Thus redistributing $\frac{q-5}{2}|S|$ into the summation over the paths in $G[S]$ we get
\begin{equation}\label{basic_step_8_main5}
3n - 3|S| - \Delta|S| - \frac{q-5}{2}|S| \leq \beta_2(S) - \sum_{P \text{ is a path in } G[S]} \left(2l(P) + (l(P)+1)\frac{q-5}{2}\right)
\end{equation}
Hence if the following inequality is satisfied
\begin{equation}\label{basic_step_9_main5}
\beta_2(S) \leq \sum_{P \text{ is a path in } G[S]} \left(2l(P) + (l(P)+1)\frac{q-5}{2}\right)
\end{equation}
Then we are done as from Inequalities \ref{basic_step_8_main5} and \ref{basic_step_9_main5} we get
\begin{equation}\label{basic_step_10_main5}
3n - 3|S| - \Delta|S| - \frac{q-5}{2}|S| \leq 0
\end{equation}
And thus
\[
|S| \geq \frac{6n}{2\Delta+q+1}
\]
The rest of the section is devoted to the proof of Inequality \ref{basic_step_9_main5}. \\
Let $T$ be the set of vertices in $V \backslash S$ that have exactly $2$ adjacent vertices in $S$.
Notice that $|T| =  \beta_2(S)$.
Given a vertex $s$ in $S$ we denote by $d_S(s)$ the degree of vertex $s$ in $G[S]$.
We shall need the following observations. \\[12pt]
\textbf{Observation 1:} If vertex $v \in T$ is adjacent to a vertex $s$ in $S$ then $d_S(s) \leq 1$. \\
\textbf{Proof:} Let $s_1$,$s_2$ be the vertices adjacent to $v$ in $S$. We consider three cases, getting a contradiction in each such case.
\begin{enumerate}
  \item Assume that $d_S(s_1) = d_S(s_2) = 2$. If $s_1$ and $s_2$ are adjacent then we remove $s_1$ from $S$ and add $v$ to $S$, thus getting a contradiction to condition $(1)$. Otherwise $s_1$ and $s_2$ are not adjacent, hence we can remove $s_1$,$s_2$ from $S$ and add $v$ to $S$, thus getting a contradiction to condition $(1)$ once again.
  \item Assume that $d_S(s_1) \leq 1$ and $d_S(s_2) = 2$. If $s_1$ and $s_2$ are on the same path in $S$ then we remove $s_2$ from $S$ and add $v$ to $S$, thus getting a contradiction to condition $(1)$. \\ Hence we can assume that $s_1$ is an endpoint of path $P_1$ in $S$ and $s_2$ is a vertex belonging to path $P_2$ in $S$, such that $P_1$ and $P_2$ are different paths in $G[S]$. If path $P_1$ is a path of length at most $3$ in $S$ then we remove $s_2$ from $S$ and add $v$ to $S$, thus getting a contradiction to condition $(1)$.
      Hence we may assume that path $P_1$ is of length $4$. Now we remove $s_1$,$s_2$ from $S$ and add $v$ to $S$, thus getting a contradiction to condition $(2)$, and we are done.
      \item Assume that $d_S(s_2) \leq 1$ and $d_S(s_1) = 2$. This case is identical to the previous
  one by symmetry.
\end{enumerate}
We conclude that $d_S(s_1) \leq 1$  and $d_S(s_2) \leq 1$. \qed \text{} \\
\noindent\textbf{Observation 2:} If a vertex $s \in S$ satisfies $d_S(s) = 1$,
then vertex $s$ has at most $q-2$ adjacent vertices in $T$. \\
\textbf{Proof:} Let $s \in S$ be a vertex which satisfies $d_S(s) = 1$. It is sufficient to prove that if vertices $v_1 \in T$ and $v_2 \in T$ are adjacent to $s$ then
 vertices $v_1$ and $v_2$ are adjacent. This implies that vertex $s$ has at most $q-2$ adjacent vertices in $T$, since graph $G$ has no cliques of size $q$. \\
We assume by contradiction that there are vertices $v_1 \in T$ and $v_2 \in T$ which adjacent to $s$ such that $v_1$ and $v_2$ are not adjacent. \\
Let $s_1 \in S$ be the second vertex adjacent to $v_1$ in $S$ (the first one being $s$).
Let $s_2 \in S$ be the second vertex adjacent to $v_2$ in $S$ (the first one being $s$).
If $s_1 = s_2$ (that is $s_1$ and $s_2$ are in fact the same vertex) then we remove $s_1,s$ from $S$ and add $v_1,v_2$ to $S$, thus getting a contradiction to condition $(1)$.
Henceforth we assume that $s_1$ and $s_2$ and different vertices. \\
Notice that by Observation 1 we have $d_S(s_1) \leq 1$ and $d_S(s_2) \leq 1$.
If $s_1$ and $s_2$ are on the same path $P$ in $G[S]$, then this path must be of length at least $1$ (as $s_1$ and $s_2$ are different vertices) and furthermore $s_1$ and $s_2$ are the two endpoints of the path (by Observation 1).
We remove vertices $s,s_1$ and add vertices $v_1,v_2$, thus getting a contradiction to condition $(1)$. \\
Finally we assume that vertex $s_1$ is in path $P_1$ and vertex $s_2$ is in path $P_2$, where $P_1$ and $P_2$ are different paths in $G[S]$.
Once again by Observation $1$ we have that vertex $s_1$ is an endpoint of path $P_1$ and vertex $s_2$ is an endpoint of path $P_2$.
First consider the case where $s$ is an endpoint of path $P_2$ (that is $s$ and $s_2$ are the two endpoints of path $P_2$).
We have the following two cases.
\begin{enumerate}
  \item Assume that path $P_1$ is of length at most $3$. We remove vertex $s$ and add vertices $v_1,v_2$, thus getting a contradiction to condition $(1)$.
  \item Assume that path $P_1$ is of length $4$. We remove vertices $s,s_1$ and add vertices $v_1,v_2$, thus getting a contradiction to condition $(1)$.
\end{enumerate}
The case of $s$ being an endpoint of path $P_1$ is handled in the same manner.
Henceforth we may assume that paths $P_1$ and $P_2$ do not contain vertex $s$.
We consider the following four cases.
\begin{enumerate}
  \item Assume that paths $P_1$ and $P_2$ are of length at most $3$. We remove $s$ from $S$ and add $v_1,v_2$ to $S$, thus getting a contradiction to condition $(1)$.
  \item Assume that paths $P_1$ and $P_2$ are of length $4$. We remove $s,s_1,s_2$ from $S$ and add $v_1,v_2$ to $S$, thus getting a contradiction to condition $(2)$.
  \item Assume that path $P_1$ is of length at most $3$ and path $P_2$ is of length $4$. We remove $s,s_2$ from $S$ and add $v_1,
  v_2$ to $S$, thus getting a contradiction to condition $(1)$.
  \item Assume that path $P_2$ is of length at most $3$ and path $P_1$ is of length $4$. We remove $s,s_1$ from $S$ and add $v_1,
  v_2$ to $S$, thus getting a contradiction to condition $(1)$.
\end{enumerate}
And thus Observation 2 follows. \qed \text{} \\
\noindent\textbf{Observation 3:} For every vertex $v \in T$ one of the following two statements holds.
\begin{enumerate}
  \item Vertex $v$ is adjacent to an endpoint of a path of length at least $3$ in $G[S]$.
  \item There is a path $P$ in $G[S]$ of length $1$ or $2$, such that $v$ is adjacent to both endpoints of $P$.
\end{enumerate}
\textbf{Proof:} By Observation $1$, vertex $v$ is adjacent to two endpoints of a single path in $G[S]$ or vertex $v$ is adjacent to endpoints of two different paths in $G[S]$.
Assume that $v$ is not adjacent to an endpoint of a path of length at least $3$ in $G[S]$.  If there is a path $P$ in $G[S]$, such that $v$ is adjacent to the two (different) endpoints of $P$ then we are done (as such path is of length $1$ or $2$). \\
Henceforth we may assume by contradiction that vertex $v$ is adjacent to vertex $s_1 \in S$ which is an endpoint of path $P_1$ and that vertex $v$ is adjacent to vertex $s_2 \in S$ which is an endpoint of path $P_2$ where $P_1$ and $P_2$ are two different paths in $G[S]$. We have three cases to consider in the following order.
\begin{enumerate}
  \item Assume that path $P_1$ is of length $0$ or path $P_2$ is of length $0$. We add vertex $v$ to $S$ thus getting a contradiction to condition $(1)$.
  \item Assume that path $P_1$ is of length $2$ or path $P_2$ is of length $2$. Let $P_1$ be a path of length $2$ without loss of generality. Remove the vertex adjacent to $s_1$ in $G[S]$ and add $v$ to $S$, thus getting a contradiction to condition $(2)$.
  \item Assume that path $P_1$ is of length $1$ and path $P_2$ is of length $1$. We add vertex $v$ to $S$, thus getting a contradiction to condition $(1)$.
\end{enumerate}
And thus Observation 3 follows. \qed \text{} \\
\noindent\textbf{Observation 4:} For any path $P$ of length $1$ in $G[S]$, there are at most $q-3$ vertices in $T$ which are adjacent to both endpoints of $P$. \\
\textbf{Proof:} Assume by contradiction that given a path $P$ of length $1$ in $G[S]$ there are at least $q-2$ vertices in $T$ which are adjacent to both endpoints of $P$.
Let $s_1,s_2$ be the two endpoints of path $P$.
As graph $G$ has no cliques of size $q$ and there are at least $q-2$ vertices in $T$ adjacent to $s_1$ and $s_2$,
there must be two vertices $v_1,v_2 \in T$ which are not adjacent such that $v_1$ is adjacent to $s_1$ and $s_2$, and $v_2$ is adjacent to $s_1$ and $s_2$.
We remove vertices $s_1,s_2$ from $S$ and add vertices $v_1,v_2$ to $S$, thus getting a contradiction to condition $(1)$.
And thus Observation 4 follows.
\qed \text{} \\
Now we are ready to prove Inequality \ref{basic_step_9_main5}. Recall that we need to prove the following.
\[
\beta_2(S) \leq \sum_{P \text{ is a path in } G[S]} \left(2l(P) + (l(P)+1)\frac{q-5}{2}\right)
\]
We will assign potential of $2l(P) + (l(P)+1)\frac{q-5}{2}$ to each path $P$ in $G[S]$.
That is the following holds:
\begin{itemize}
  \item A path of length at least $3$ has a potential of at least $2q-4$.
  \item A path of length $2$ has a potential of $\frac{3q-7}{2}$.
  \item A path of length $1$ has a potential of $q-3$.
\end{itemize}
Now we shall show how to redistribute this potential as to give to each vertex in $T$ at least one unit of potential. \\
Notice that by Observation 1 vertices in $T$ can be adjacent only to endpoints of paths in $G[S]$.
If path $P$ is of length at least $3$ then its potential is at least $2q-4 = 2(q-2)$ and furthermore by Observation 2 there are at most $2(q-2)$ vertices
in $T$ which are adjacent to an endpoint of path $P$ and thus we may give each such adjacent vertex in $T$ a potential of $1$. \\
By Observation 3 every vertex in $T$ which is not adjacent to a path of length at least $3$ is adjacent to both endpoints of some path $P$ in $G[S]$,
where $P$ is of length $1$ or $2$. \\
By Observation 2, given a path $P$ of length $2$ in $G[S]$, there are at most $q-2$ vertices in $T$ that are adjacent to both endpoints of $P$.
Hence each such path $P$ can contribute a potential of at least $\frac{3q-7}{2(q-2)} \geq 1$ to each of the vertices of $T$ that are adjacent to both endpoints of $P$. \\
Finally by Observation 4, given a path $P$ of length $1$ in $G[S]$, there are at most $q-3$ vertices in $T$ that are adjacent to both endpoints of $P$.
Hence each such path $P$ can contribute a potential of $\frac{q-3}{q-3}=1$ to each of the vertices of $T$ that are adjacent to both endpoints of $P$.
We showed that each vertex in $T$ gets a potential of at least $1$ and the proof follows.
\EPF

\section{Proof of Theorem \ref{beautiful_theorem1}}
We shall prove the following theorem. \\
Let $G$ be a graph of order $n$, maximum degree $\Delta>0$ and maximum clique size $\omega$. Then
\[
a(G) \geq \frac{6n}{2\Delta + \omega +2}.
\]
We have $3$ cases:
\begin{itemize}
  \item If $\omega = 2$ then the theorem follows from Corollary \ref{our_bound_triangle_free_cor}.
  \item If $\omega = 3$ then the theorem follows from Theorem \ref{four_clique_main_thm}.
  \item If $\omega \geq 4$ then the theorem follows from Theorem \ref{five_theorem_main}.
\end{itemize}
And thus Theorem \ref{beautiful_theorem1} is proven.

\appendix
\section{Appendix A}\label{app1sect}
Recall that a linear $k$-forest is a forest consisting of paths of length at most $k$, and that $a_k(G)$ denote the maximum size of an induced linear $k$-forest
in $G$.
The following bound was first proven in \cite{hopkins1}
(It is a straightforward corollary of \cite{lovlemma1}).
\BTHM
Let $G=(V,E)$ be a graph of order $n$ and maximum degree $\Delta$ where $\Delta$ is odd. Then
  \[
  a_1(G) \geq \frac{2n}{\Delta+1}.
  \]
\ETHM
We will prove the following theorem in this appendix.
\BTHM\label{appendix_tiny_thm_1}
Let $G=(V,E)$ be a graph of order $n$ and maximum degree $\Delta>0$. Then
  \[
  a_3(G) \geq \frac{2n}{\Delta+1}.
  \]
\ETHM
We will start with the following lemma.
\BL\label{appendix_tiny_lemma_1}
Let $G=(V,E)$ be a $\Delta$-regular graph of order $n$ where $\Delta>0$. Then
  \[
  a_3(G) \geq \frac{2n}{\Delta+1}.
  \]
\EL
\BPF
Given a set $S$ of vertices of $G$, Let $|S|$ denote the number of vertices in $S$, $G[S]$ denote the subgraph of $G$ induced by the vertices of $S$, and
$e(S)$ denote the number of edges in $G[S]$. \\
Choose an induced linear $3$-forest $S$ in graph $G$ such that the following conditions are satisfied.
\begin{enumerate}[label=(\arabic*), leftmargin=1.5cm]
  \item $|S|$ is maximized.
  \item Subject to $(1)$, $e(S)$ is minimized.
  \item Subject to $(2)$, the number of vertices of degree $1$ in $G[S]$ is maximized.
\end{enumerate}
Let $\beta_i(S)$ denote the number of vertices in $V \backslash S$ with exactly $i$ adjacent vertices in $S$.
Notice that
\begin{equation}\label{basic_step_1_appendix}
  n - |S| = \sum_{i=0}^{|S|} \beta_i(S).
\end{equation}
Given a vertex $s \in S$, let $d^o(s)$ denote the number of vertices in $V \backslash S$ that are adjacent to $s$.
Notice that
\begin{equation}\label{basic_step_2_appendix}
 \sum_{s \in S} d^o(s) = \sum_{i=0}^{|S|} i \beta_i(S).
\end{equation}
Multiplying Equation \ref{basic_step_1_appendix} by $2$ and subtracting Equation \ref{basic_step_2_appendix} we obtain the following.
\begin{equation}\label{basic_step_3_appendix}
2n - 2|S| - \sum_{s \in S} d^o(s) =  \sum_{i=0}^{|S|} (2 - i) \beta_i(S).
\end{equation}
By condition $(1)$ we have $\beta_0(S) = 0$.  Now we shall prove that $\beta_1(S) = 0$. \\
Assume by contradiction that there is a vertex $v$ in $V \backslash S$ with exactly one adjacent vertex $s$ in $S$.
If $s$ is a vertex of degree at least $1$ in $G[S]$ then we remove vertex $s$ from $S$ and add vertex $v$ to $S$, thus getting a contradiction to condition $(2)$.
Hence $s$ is a vertex of degree $0$ in $G[S]$. Thus we can add $v$ to $S$ and get a contradiction to condition $(1)$. \\
We have shown that $\beta_0(S) = 0$ and $\beta_1(S) = 0$. Thus from equality \ref{basic_step_3_appendix} it follows that
\begin{equation}\label{basic_step_4_appendix}
2n =  2|S| + \sum_{s \in S} d^o(s) - \sum_{i=3}^{|S|} (i - 2) \beta_i(S).
\end{equation}
Let $S_0$ be the set of vertices of degree $0$ in $G[S]$ and set $S_1 = S \backslash S_0$.
As the degree of each vertex of $S_1$ in $G[S]$ is positive we have the following inequality.
\begin{equation}\label{basic_step_5_inequality_appendix}
\sum_{s \in S} d^o(s) \leq \Delta |S| - |S_1|.
\end{equation}
Thus
\begin{equation}\label{basic_step_6_appendix}
2n \leq  2|S| + \Delta |S| - |S_1| - \sum_{i=3}^{|S|} (i - 2) \beta_i(S).
\end{equation}
Now notice that it follows from \ref{basic_step_6_appendix} that if $\sum_{i=3}^{|S|} (i - 2) \beta_i(S) \geq |S_0|$
then $|S| \geq \frac{2n}{\Delta+1}$.
Thus all that remains is to show that (for $\Delta >0$)
\begin{equation}\label{basic_step_7_appendix}
\sum_{i=3}^{|S|} (i - 2) \beta_i(S) \geq |S_0|.
\end{equation}
In fact we shall show
\begin{equation}\label{basic_step_8_appendix}
\sum_{i=3}^{|S|} \beta_i(S) \geq |S_0|.
\end{equation}
Notice that Inequality \ref{basic_step_7_appendix} follows from \ref{basic_step_8_appendix}. \\
Let $T$ be the set of vertices in $V \backslash S$ that have at least $3$ adjacent vertices in $S$.
Notice that $|T| = \sum_{i=3}^{|S|} \beta_i(S)$.
We shall need the following observation. \\[12pt]
\textbf{Observation I:} if vertex $v \in V \backslash S$ is adjacent to a vertex in $S_0$ then  $v \in T$. \\
\textbf{Proof:} Notice that each vertex $v \in V \backslash S$ has at least two adjacent vertices in $S$ as $\beta_0(S) = 0$ and $\beta_1(S) = 0$.
Assume by contradiction that there is a vertex $v$ in $V \backslash S$ with exactly two adjacent vertex $s_0$ and $s_1$ in $S$, such that
vertex $s_0$ is in $S_0$.
If $s_1$ is a vertex of degree $2$ in $G[S]$ then we remove vertex $s_1$ from $S$ and add vertex $v$ to $S$, thus getting a contradiction
to condition $(2)$. Hence $s_1$ is an endpoint of some path $P$ in $S$. If path $P$ is of length at most $1$ then we can add vertex $v$ to $S$ and get a contradiction to condition $(1)$.
Thus path $P$ is of length at least $2$. Now we remove vertex $s_1$ from $S$ and add vertex $v$ to $S$, thus getting a contradiction to condition $(3)$.
And thus observation I follows. \qed \text{} \\
By the regularity of $G$ and Observation I we have that each vertex $s \in S_0$ has exactly $\Delta$ adjacent vertices in $T$.
Hence as we have $\Delta |S_0|$ edges between $S_0$ and $T$, we conclude that $|T| \geq |S_0|$ and thus $\sum_{i=3}^{|S|} \beta_i(S) = |T| \geq |S_0|$ and we are done.
\EPF
We shall prove now that Theorem \ref{appendix_tiny_thm_1} follows from Lemma \ref{appendix_tiny_lemma_1}. \\
Observe that, if $G$ is a graph with maximum degree $\Delta>0$, then we can create a $\Delta$-regular
graph by taking copies $H_1,H_2, \dots , H_r$ of $G$ and joining some pairs of vertices from different copies so
as to make the resulting graph $G'$ a $\Delta$-regular graph. Applying Lemma \ref{appendix_tiny_lemma_1} to graph $G'$
we get by the pigeonhole principle that for some $1 \leq i \leq r$ we have $a_3(H_i) \geq \frac{2n}{\Delta+1}$ and thus we are done.

\section{Appendix B}\label{app2sect}
We shall prove in this appendix Observation $4$ of Section \ref{four_clique_section}.
That is we shall prove that a tree $T$ in $G[S]$ on at most $7$ vertices has at most $|T|-1$ adjacent vertices in $B_2$.
Recall that tree $T$ can have at most $|T|$ adjacent vertices in $B_2$ (Observation $3$ in Section \ref{four_clique_section}).
Assume by contradiction that tree $T$ has of exactly $|T|$ adjacent vertices in $B_2$. Let $A \subseteq B_2$ be the set of $|T|$ vertices that are adjacent to tree $T$.
Recall that each vertex in $B_2$ is either adjacent to two vertices in $T$ or not adjacent to any vertex in $T$, for otherwise we get a contradiction to condition $(1)$
 in Section \ref{four_clique_section}.
Now notice that each vertex in $A$ has exactly two adjacent vertices in $T$ and every vertex in $T$ has exactly two adjacent vertices in $A$
(this follows from Observation $2$ in Section \ref{four_clique_section}). \\
Recall that we have chosen in Section \ref{four_clique_section} an induced forest $S$ in graph $G$ such that the following conditions are satisfied.
\begin{enumerate}[label=(\arabic*), leftmargin=1.5cm]
  \item $|S|$ is maximized.
  \item Subject to $(1)$, $e(S)$ is maximized.
  \item Subject to $(2)$, the number of vertices of degree $1$ in $G[S]$ is maximized.
  \item  Subject to $(3)$, we maximize the following sum.
  \[
  \sum_{T \text{ is a tree in } G[S]} \Delta(T).
  \]
  \item Subject to $(4)$, we minimize the following sum.
  \[
  \sum_{T \text{ is a tree in } G[S]} P(T).
  \]
\end{enumerate}
Recall that $e(S)$ denotes the number of edges in $G[S]$. Given an induced subgraph $T$ of $G[S]$ we denote by $\Delta(T)$ the maximum degree of $T$.
We denote by $D(T)$ the diameter of $T$ (that is the great distance between any pair of vertices in $T$).
Finally we denote by $P(T)$ the number of paths in $T$ of length $D(T)$. \\
We shall need a few claims.
\\ \text{} \\
\noindent\textbf{Claim 1:} If vertices $s_1$ and $s_2$ in $T$ are adjacent then there is at most one vertex $v \in A$ such that $v$ is adjacent both to $s_1$ and $s_2$. \\
\textbf{Proof:} Assume that there are vertices $v_1,v_2$  in $A$ that are adjacent to both vertices $s_1,s_2$ in $T$.
Vertices $v_1$ and $v_2$ can not be adjacent as graph $G$ has no cliques of size $4$, hence we can add vertices $v_1,v_2$ to $S$ and remove vertex $s_1$ from $S$, thus getting a contradiction
to condition $(1)$. \qed \text{} \\
\noindent\textbf{Claim 2:} Let $s$ be a vertex in $T$ and let $v_1,v_2$ be the vertices adjacent to $s$ in $A$. Let $s_1$ be the second neighbor of $v_1$ in $T$ and let $s_2$ be second neighbor of $v_2$ in $T$. Remove vertex $s$ from tree $T$ and denote the resulting forest by $T'$. Then vertices $s_1$ and $s_2$ belong to the same connected component in $T'$. \\
\textbf{Proof:} If vertices $s_1$ and $s_2$ belong to different connected components of $T'$ then we can add vertices $v_1,v_2$ to $S$ and remove vertex $s$ from $S$ thus getting a contradiction to condition $(1)$.
 \qed \text{} \\
 \noindent\textbf{Claim 3:} Let $s$ be a vertex in $T$ and let $v_1,v_2$ be the vertices adjacent to $s$ in $A$. Let $s_1$ be the second neighbor of $v_1$ in $T$ and let $s_2$ be second neighbor of $v_2$ in $T$. Then vertex $s$ can not be adjacent to both vertices $s_1$ and $s_2$. \\
\textbf{Proof:} Assume by contradiction that $s$ is adjacent to $s_1$ and $s_2$.
By Claim $1$ we have that $s_1 \neq s_2$, but then we get a contradiction to Claim $2$ for vertex $s$.
 \qed \text{} \\
\noindent\textbf{Claim 4:} Let $s$ be a leaf vertex in $T$ (that is $d_T(s)=1$). If vertex $s$ is adjacent to a vertex $s_1$ in $T$ such that $d_T(s_1) = 2$ then
for any vertex $v \in A$ that is adjacent to $s$ , the second neighbor of $v$ in $T$ must be a leaf vertex too. \\
\textbf{Proof:} Let $s_2$ be the second neighbor of $v$ in the tree $T$ (the first neighbor is $s$). if $d_T(s_2) \geq 2$ and $s_2 \neq s_1$ then we can remove
vertex $s$ from $S$ and add vertex $v$ to $S$, thus getting a contradiction to condition $(3)$. \\
If $s_2 = s_1$ then by Claim $1$ vertex $v$ must be the only vertex in $A$ that is adjacent both to $s$ and $s_2$. Hence vertex $s_2$ has an adjacent vertex $v_2$ in $A$ such that the
second neighbor of $v_2$ in $T$ is a vertex different from $s$, but that is a contradiction to Claim $2$ (for vertex $s_2$) and thus we are done. \qed \text{} \\
Now we shall do a case analysis on all non-isomorphic trees of at most $7$ vertices. \\
\textbf{Case $1$: Tree $T$ is an isolated vertex}. We get a contradiction to condition $(1)$ \\
\noindent\textbf{Case $2$: Tree $T$ is a star (on any number of vertices)}. \\
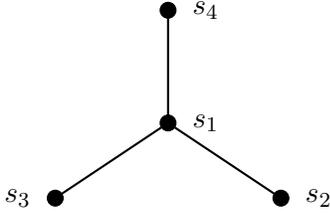
\begin{tikzpicture}
\draw[fill=black] (0,0) circle (3pt);
\draw[fill=black] (3,0) circle (3pt);
\draw[fill=black] (1.5,1) circle (3pt);
\draw[fill=black] (1.5,2.5) circle (3pt);
\node at (-0.5,0) {$s_3$};
\node at (3.5,0) {$s_2$};
\node at (2,1) {$s_1$};
\node at (2,2.5) {$s_4$};
\draw[thick] (1.5,1) -- (1.5,2.5);
\draw[thick] (1.5,1) -- (3,0);
\draw[thick] (1.5,1) -- (0,0);
\end{tikzpicture} \\
Let $s_1$ be the center of the star. We get a contradiction to Claim $3$ for vertex $s_1$. \\
\noindent\textbf{Case 4.1: Tree $T$ is the following tree}. \\
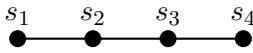
\begin{tikzpicture}
\draw[fill=black] (0,0) circle (3pt);
\draw[fill=black] (1,0) circle (3pt);
\draw[fill=black] (2,0) circle (3pt);
\draw[fill=black] (3,0) circle (3pt);
\node at (0,0.3) {$s_1$};
\node at (1,0.3) {$s_2$};
\node at (2,0.3) {$s_3$};
\node at (3,0.3) {$s_4$};
\draw[thick] (0,0) -- (1,0) -- (2,0) -- (3,0);
\end{tikzpicture} \\
By Claim $4$ there are vertices $v_1,v_2 \in A$ that are adjacent to both vertices $s_1,s_4 \in T$.
Hence we get a contradiction to Claim $1$ for vertices $s_2,s_3$. \\
\noindent\textbf{Case 5.1: Tree $T$ is the following tree}. \\
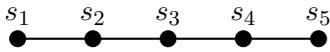
\begin{tikzpicture}
\draw[fill=black] (0,0) circle (3pt);
\draw[fill=black] (1,0) circle (3pt);
\draw[fill=black] (2,0) circle (3pt);
\draw[fill=black] (3,0) circle (3pt);
\draw[fill=black] (4,0) circle (3pt);
\node at (0,0.3) {$s_1$};
\node at (1,0.3) {$s_2$};
\node at (2,0.3) {$s_3$};
\node at (3,0.3) {$s_4$};
\node at (4,0.3) {$s_5$};
\draw[thick] (0,0) -- (1,0) -- (2,0) -- (3,0) -- (4,0);
\end{tikzpicture} \\
By Claim $4$ there are vertices $v_1,v_2 \in A$ that are adjacent to both vertices $s_1,s_5 \in T$.
Hence we get a contradiction to Claim $3$ for the vertex $s_3$. \\
\noindent\textbf{Case 5.2: Tree $T$ is the following tree}. \\
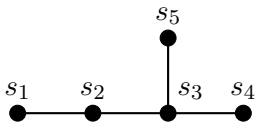
\begin{tikzpicture}
\draw[fill=black] (0,0) circle (3pt);
\draw[fill=black] (1,0) circle (3pt);
\draw[fill=black] (2,0) circle (3pt);
\draw[fill=black] (3,0) circle (3pt);
\draw[fill=black] (2,1) circle (3pt);
\node at (0,0.3) {$s_1$};
\node at (1,0.3) {$s_2$};
\node at (2.3,0.3) {$s_3$};
\node at (3,0.3) {$s_4$};
\node at (2,1.3) {$s_5$};
\draw[thick] (0,0) -- (1,0) -- (2,0) -- (3,0);
\draw[thick] (2,0) -- (2,1);
\end{tikzpicture} \\
By Claim $4$ there is no vertex in $A$ that is adjacent to both $s_3$ and $s_1$.
Hence we get a contradiction to Claim $3$ for the vertex $s_3$. \\
\noindent\textbf{Case 6.1: Tree $T$ is the following tree}. \\
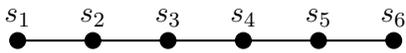
\begin{tikzpicture}
\draw[fill=black] (0,0) circle (3pt);
\draw[fill=black] (1,0) circle (3pt);
\draw[fill=black] (2,0) circle (3pt);
\draw[fill=black] (3,0) circle (3pt);
\draw[fill=black] (4,0) circle (3pt);
\draw[fill=black] (5,0) circle (3pt);
\node at (0,0.3) {$s_1$};
\node at (1,0.3) {$s_2$};
\node at (2,0.3) {$s_3$};
\node at (3,0.3) {$s_4$};
\node at (4,0.3) {$s_5$};
\node at (5,0.3) {$s_6$};

\draw[thick] (0,0) -- (1,0) -- (2,0) -- (3,0) -- (4,0) -- (5,0);
\end{tikzpicture} \\
By Claim $4$ there are vertices $v_1,v_2 \in A$ that are adjacent to both vertices $s_1,s_6 \in T$.
Since $G$ has no cliques of size $4$ there are two non-adjacent vertices $v_3,v_4 \in A \backslash \{v_1,v_2\}$.
If there is a vertex $s$ in $T$ that is adjacent to both $v_3$ and $v_4$ then we remove vertex $s$ from $S$ and add vertices $v_3,v_4$ to $S$,
thus getting a contradiction to condition $(1)$.
Assume w.l.o.g. that $v_3$ is adjacent to $s_4$. Now we have the following cases.
\begin{enumerate}
  \item If  $v_3$ is adjacent to $s_4,s_5$ then by Claim $1$ there is a vertex $u \in A$ that is adjacent to $s_4$ and $s_i$ for some $2 \leq i \leq 3$.
Hence we get a contradiction to Claim $2$ for the vertex $s_4$.
  \item If  $v_3$ is adjacent to $s_4,s_3$ then $v_4$ adjacent to $s_2,s_5$. We remove vertex $s_3$ from $S$ and add vertices $v_3,v_4$ to $S$, thus getting a contradiction to condition $(1)$.
  \item It $v_3$ is adjacent to $s_2,s_4$ then $v_4$ is adjacent to $s_3,s_5$. We remove vertex $s_3$ from $S$ and add vertices $v_3,v_4$ to $S$, thus getting a contradiction to condition $(1)$.
\end{enumerate}
\noindent\textbf{Case 6.2: Tree $T$ is the following tree}. \\
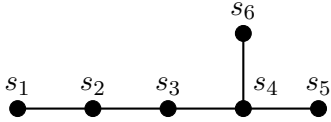
\begin{tikzpicture}
\draw[fill=black] (0,0) circle (3pt);
\draw[fill=black] (1,0) circle (3pt);
\draw[fill=black] (2,0) circle (3pt);
\draw[fill=black] (3,0) circle (3pt);
\draw[fill=black] (4,0) circle (3pt);
\draw[fill=black] (3,1) circle (3pt);
\node at (0,0.3) {$s_1$};
\node at (1,0.3) {$s_2$};
\node at (2,0.3) {$s_3$};
\node at (3.3,0.3) {$s_4$};
\node at (4,0.3) {$s_5$};
\node at (3,1.3) {$s_6$};

\draw[thick] (0,0) -- (1,0) -- (2,0) -- (3,0) -- (4,0);
\draw[thick] (3,0) -- (3,1);
\end{tikzpicture} \\
Let vertex $v \in A$ be adjacent to vertex $s_1 \in T$.
By Claim $4$ the second neighbor of $v$ in $T$ is either $s_5$ or $s_6$.
We remove vertex $s_1$ from $S$ and add vertex $v$ to $S$, thus getting a contradiction to condition $(5)$. \\
\noindent\textbf{Case 6.3: Tree $T$ is the following tree}. \\
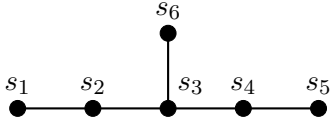
\begin{tikzpicture}
\draw[fill=black] (0,0) circle (3pt);
\draw[fill=black] (1,0) circle (3pt);
\draw[fill=black] (2,0) circle (3pt);
\draw[fill=black] (3,0) circle (3pt);
\draw[fill=black] (4,0) circle (3pt);
\draw[fill=black] (2,1) circle (3pt);
\node at (0,0.3) {$s_1$};
\node at (1,0.3) {$s_2$};
\node at (2.3,0.3) {$s_3$};
\node at (3,0.3) {$s_4$};
\node at (4,0.3) {$s_5$};
\node at (2,1.3) {$s_6$};

\draw[thick] (0,0) -- (1,0) -- (2,0) -- (3,0) -- (4,0);
\draw[thick] (2,0) -- (2,1);
\end{tikzpicture} \\
Let vertex $v \in A$ be adjacent to vertex $s_3 \in T$.
By Claim $4$ the second neighbor of $v$ in $T$ can not be $s_1$ or $s_5$.
Hence we get a contradiction to Claim $3$ for vertex $s_3$. \\
\noindent\textbf{Case 6.4: Tree $T$ is the following tree}. \\
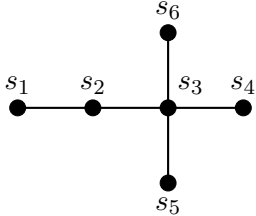
\begin{tikzpicture}
\draw[fill=black] (0,0) circle (3pt);
\draw[fill=black] (1,0) circle (3pt);
\draw[fill=black] (2,0) circle (3pt);
\draw[fill=black] (3,0) circle (3pt);
\draw[fill=black] (2,-1) circle (3pt);
\draw[fill=black] (2,1) circle (3pt);
\node at (0,0.3) {$s_1$};
\node at (1,0.3) {$s_2$};
\node at (2.3,0.3) {$s_3$};
\node at (3,0.3) {$s_4$};
\node at (2,-1.3) {$s_5$};
\node at (2,1.3) {$s_6$};

\draw[thick] (0,0) -- (1,0) -- (2,0) -- (3,0);
\draw[thick] (2,0) -- (2,1);
\draw[thick] (2,0) -- (2,-1);
\end{tikzpicture} \\
Let vertex $v \in A$ be adjacent to vertex $s_3 \in T$.
By Claim $4$ the second neighbor of $v$ in $T$ can not be $s_1$.
Hence we get a contradiction to Claim $3$ for vertex $s_3$. \\
\noindent\textbf{Case 6.5: Tree $T$ is the following tree}. \\
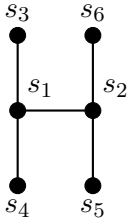
\begin{tikzpicture}
\draw[fill=black] (0,0) circle (3pt);
\draw[fill=black] (1,0) circle (3pt);
\draw[fill=black] (0,1) circle (3pt);
\draw[fill=black] (0,-1) circle (3pt);
\draw[fill=black] (1,-1) circle (3pt);
\draw[fill=black] (1,1) circle (3pt);
\node at (0.3,0.3) {$s_1$};
\node at (1.3,0.3) {$s_2$};
\node at (0,1.3) {$s_3$};
\node at (0,-1.3) {$s_4$};
\node at (1,-1.3) {$s_5$};
\node at (1,1.3) {$s_6$};

\draw[thick] (0,0) -- (1,0);
\draw[thick] (0,-1) -- (0,0) -- (0,1);
\draw[thick] (1,-1) -- (1,0) -- (1,1);
\end{tikzpicture} \\
Let vertex $v \in A$ be adjacent to vertex $s_1 \in T$.
The second neighbor of $v$ in $T$ can not be $s_5$. As if $v$ is adjacent to $s_5$ we can remove $s_5$ from $S$ and add $v$ to $S$, thus getting a contradiction
to condition $(4)$. In the same manner we can show that the second neighbor of $v$ in $T$ can not be $s_6$.
Hence we get a contradiction to Claim $3$ for vertex $s_1$. \\
\noindent\textbf{Case 7.1: Tree $T$ is the following tree}. \\
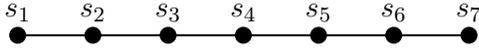
\begin{tikzpicture}
\draw[fill=black] (0,0) circle (3pt);
\draw[fill=black] (1,0) circle (3pt);
\draw[fill=black] (2,0) circle (3pt);
\draw[fill=black] (3,0) circle (3pt);
\draw[fill=black] (4,0) circle (3pt);
\draw[fill=black] (5,0) circle (3pt);
\draw[fill=black] (6,0) circle (3pt);
\node at (0,0.3) {$s_1$};
\node at (1,0.3) {$s_2$};
\node at (2,0.3) {$s_3$};
\node at (3,0.3) {$s_4$};
\node at (4,0.3) {$s_5$};
\node at (5,0.3) {$s_6$};
\node at (6,0.3) {$s_7$};

\draw[thick] (0,0) -- (1,0) -- (2,0) -- (3,0) -- (4,0) -- (5,0) -- (6,0);
\end{tikzpicture} \\
By Claim $4$ there are vertices $v_1,v_2 \in A$ that are adjacent to both vertices $s_1,s_7 \in T$.
Since $G$ has no cliques of size $4$ there are two non-adjacent vertices $v_3,v_4 \in A \backslash \{v_1,v_2\}$.
If there is a vertex $s$ in $T$ that is adjacent to both $v_3$ and $v_4$ then we remove vertex $s$ from $S$ and add vertices $v_3,v_4$ to $S$,
thus getting a contradiction to condition $(1)$. \\
Notice that there is no vertex $u \in A$ such that $u$ is adjacent to $s_2,s_3$ since if 
 $u$ is adjacent to $s_2,s_3$ then by Claim $1$ there is a vertex $u' \in A$ that is adjacent to $s_3$ and $s_i$ for some $4 \leq i \leq 6$.
Hence we get a contradiction to Claim $2$ for the vertex $s_3$. By the same logic there is no vertex $u \in A$ such that $u$ is adjacent to $s_5,s_6$  \\
If one of the vertices $v_3,v_4$ is adjacent to $s_4$ (assume w.l.o.g. that it is $v_3$) then we may assume by symmetry that one of the following cases occurs.
\begin{enumerate}
  \item Vertex $v_3$ is adjacent to vertices $s_3,s_4$. In this case we may assume that vertex $v_4$ is adjacent to vertices $s_2,s_5$ or 
 vertex $v_4$ is adjacent to vertices $s_2,s_6$. In both cases we remove vertex $s_4$ from $S$ and add vertices $v_3,v_4$ to $S$, thus getting a contradiction to condition $(1)$.
  \item Vertex $v_3$ is adjacent to vertices $s_2,s_4$. In this case we may assume that vertex $v_4$ is adjacent to vertices $s_3,s_5$ or vertex $v_4$ is adjacent to vertices $s_3,s_6$. In both cases we remove vertex $s_4$ from $S$ and add vertices $v_3,v_4$ to $S$, thus getting a contradiction to condition $(1)$.
\end{enumerate}
Hence we may assume that vertices  $v_3,v_4$ are not adjacent to vertex $s_4$. Thus we may assume that one of the two following cases occurs.
\begin{itemize}
  \item Vertex $v_3$ is adjacent to vertices $s_2,s_5$ and vertex $v_4$ is adjacent to vertices $s_3,s_6$.
  \item Vertex $v_3$ is adjacent to vertices $s_2,s_6$ and vertex $v_4$ is adjacent to vertices $s_3,s_5$.
\end{itemize}
In both cases we remove vertex $s_3$ from $S$ and add vertex $v_3,v_4$ to $S$, thus getting a contradiction to condition $(1)$. \\
\noindent\textbf{Case 7.2: Tree $T$ is the following tree}. \\
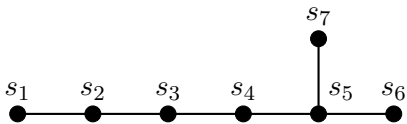
\begin{tikzpicture}
\draw[fill=black] (0,0) circle (3pt);
\draw[fill=black] (1,0) circle (3pt);
\draw[fill=black] (2,0) circle (3pt);
\draw[fill=black] (3,0) circle (3pt);
\draw[fill=black] (4,0) circle (3pt);
\draw[fill=black] (5,0) circle (3pt);
\draw[fill=black] (4,1) circle (3pt);
\node at (0,0.3) {$s_1$};
\node at (1,0.3) {$s_2$};
\node at (2,0.3) {$s_3$};
\node at (3,0.3) {$s_4$};
\node at (4.3,0.3) {$s_5$};
\node at (5,0.3) {$s_6$};
\node at (4,1.3) {$s_7$};

\draw[thick] (0,0) -- (1,0) -- (2,0) -- (3,0) -- (4,0) -- (5,0);
\draw[thick] (4,0) -- (4,1);
\end{tikzpicture} \\
Let vertex $v \in A$ be adjacent to vertex $s_1 \in T$.
By Claim $4$ the second neighbor of $v$ in $T$ is either $s_6$ or $s_7$.
We remove vertex $s_1$ from $S$ and add vertex $v$ to $S$, thus getting a contradiction to condition $(5)$. \\
\noindent\textbf{Case 7.3: Tree $T$ is the following tree}. \\
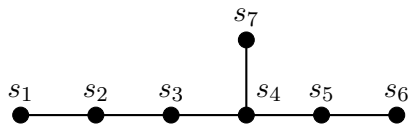
\begin{tikzpicture}
\draw[fill=black] (0,0) circle (3pt);
\draw[fill=black] (1,0) circle (3pt);
\draw[fill=black] (2,0) circle (3pt);
\draw[fill=black] (3,0) circle (3pt);
\draw[fill=black] (4,0) circle (3pt);
\draw[fill=black] (5,0) circle (3pt);
\draw[fill=black] (3,1) circle (3pt);
\node at (0,0.3) {$s_1$};
\node at (1,0.3) {$s_2$};
\node at (2,0.3) {$s_3$};
\node at (3.3,0.3) {$s_4$};
\node at (4,0.3) {$s_5$};
\node at (5,0.3) {$s_6$};
\node at (3,1.3) {$s_7$};

\draw[thick] (0,0) -- (1,0) -- (2,0) -- (3,0) -- (4,0) -- (5,0);
\draw[thick] (3,0) -- (3,1);
\end{tikzpicture} \\
By Claim $4$ one of the following two subcases occurs.
\begin{itemize}
  \item Subcase $1$: there are vertices $v_1,v_2 \in A$ that are adjacent to both vertices $s_1,s_6 \in T$.
  \item Subcase $2$: there are vertices $v_1,v_2,v_3 \in A$ such that vertex $v_1$ is adjacent to $s_1,s_6$, vertex $v_2$ is adjacent to $s_6,s_7$ and
   vertex $v_3$ is adjacent to $s_1,s_7$.
\end{itemize}
The analysis of the Subcase $2$ is identical to the analysis of Case $6.1$
(since the vertices of $A \backslash \{v_1,v_2,v_3\}$ are adjacent to the vertices of a path of length $3$ in $T$ in that case).
Hence we may assume that Subcase $1$ occurs, that is there are vertices $v_1,v_2 \in A$ that are adjacent to both vertices $s_1,s_6 \in T$.
Since $G$ has no cliques of size $4$ there are two non-adjacent vertices $v_3,v_4 \in A \backslash \{v_1,v_2\}$.
If there is a vertex $s$ in $T$ that is adjacent to both $v_3$ and $v_4$ then we remove vertex $s$ from $S$ and add vertices $v_3,v_4$ to $S$,
thus getting a contradiction to condition $(1)$. \\
Notice that there is no vertex $u \in A$ such that $u$ is adjacent to $s_2,s_3$ since if
 $u$ is adjacent to $s_2,s_3$ then by Claim $1$ there is a vertex $u' \in A$ that is adjacent to $s_3$ and $s_i$ for some $i \geq  4$.
Hence we get a contradiction to Claim $2$ for the vertex $s_3$. \\
If one of the vertices $v_3,v_4$ is adjacent to $s_4$ then we remove vertex $s_4$ from $S$ and add vertices $v_3,v_4$ to $S$, thus getting a contradiction to condition $(1)$.
Otherwise we may assume that vertex $v_3$ is adjacent to vertices $s_2,s_7$ and vertex $v_4$ is adjacent to vertices $s_3,s_5$.
We remove vertex $s_3$ from $S$ and add vertices $v_3,v_4$ to $S$, thus getting a contradiction to condition $(1)$. \\
\noindent\textbf{Case 7.4: Tree $T$ is the following tree}. \\
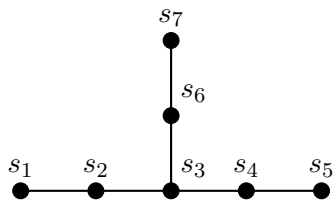
\begin{tikzpicture}
\draw[fill=black] (0,0) circle (3pt);
\draw[fill=black] (1,0) circle (3pt);
\draw[fill=black] (2,0) circle (3pt);
\draw[fill=black] (3,0) circle (3pt);
\draw[fill=black] (4,0) circle (3pt);
\draw[fill=black] (2,1) circle (3pt);
\draw[fill=black] (2,2) circle (3pt);
\node at (0,0.3) {$s_1$};
\node at (1,0.3) {$s_2$};
\node at (2.3,0.3) {$s_3$};
\node at (3,0.3) {$s_4$};
\node at (4,0.3) {$s_5$};
\node at (2.3,1.3) {$s_6$};
\node at (2,2.3) {$s_7$};
\draw[thick] (0,0) -- (1,0) -- (2,0) -- (3,0) -- (4,0);
\draw[thick] (2,0) -- (2,1) -- (2,2);
\end{tikzpicture} \\
Let vertex $v \in A$ be adjacent to vertex $s_3 \in T$.
By Claim $4$ the second neighbor of $v$ in $T$ can not be $s_1$ or $s_5$ or $s_7$.
Hence we get a contradiction to Claim $3$ for vertex $s_3$. \\
\noindent\textbf{Case 7.5: Tree $T$ is the following tree}. \\
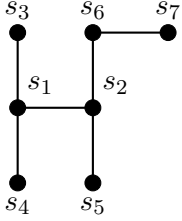
\begin{tikzpicture}
\draw[fill=black] (0,0) circle (3pt);
\draw[fill=black] (1,0) circle (3pt);
\draw[fill=black] (0,1) circle (3pt);
\draw[fill=black] (0,-1) circle (3pt);
\draw[fill=black] (1,-1) circle (3pt);
\draw[fill=black] (1,1) circle (3pt);
\draw[fill=black] (2,1) circle (3pt);
\node at (0.3,0.3) {$s_1$};
\node at (1.3,0.3) {$s_2$};
\node at (0,1.3) {$s_3$};
\node at (0,-1.3) {$s_4$};
\node at (1,-1.3) {$s_5$};
\node at (1,1.3) {$s_6$};
\node at (2,1.3) {$s_7$};

\draw[thick] (0,0) -- (1,0);
\draw[thick] (0,-1) -- (0,0) -- (0,1);
\draw[thick] (1,-1) -- (1,0) -- (1,1);
\draw[thick] (1,1) -- (2,1);
\end{tikzpicture} \\
Let vertex $v \in A$ be adjacent to vertex $s_2 \in T$.
By Claim $4$ the second neighbor of $v$ in $T$ can not be $s_7$.
Furthermore the second neighbor of $v$ in $T$ can not be $s_3$. As if $v$ is adjacent to $s_3$ we can remove $s_3$ from $S$ and add $v$ to $S$, thus getting a contradiction
to condition $(4)$. In the same manner we can show that the second neighbor of $v$ in $T$ can not be $s_4$.
Hence we get a contradiction to Claim $3$ for vertex $s_2$. \\
\noindent\textbf{Case 7.6: Tree $T$ is the following tree}. \\
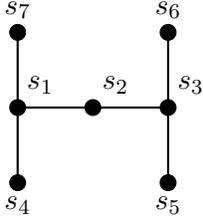
\begin{tikzpicture}
\draw[fill=black] (0,0) circle (3pt);
\draw[fill=black] (1,0) circle (3pt);
\draw[fill=black] (0,1) circle (3pt);
\draw[fill=black] (0,-1) circle (3pt);
\draw[fill=black] (2,-1) circle (3pt);
\draw[fill=black] (2,1) circle (3pt);
\draw[fill=black] (2,0) circle (3pt);
\node at (0.3,0.3) {$s_1$};
\node at (1.3,0.3) {$s_2$};
\node at (0,1.3) {$s_7$};
\node at (0,-1.3) {$s_4$};
\node at (2,-1.3) {$s_5$};
\node at (2,1.3) {$s_6$};
\node at (2.3,0.3) {$s_3$};

\draw[thick] (0,0) -- (1,0) -- (2,0);
\draw[thick] (0,-1) -- (0,0) -- (0,1);
\draw[thick] (2,-1) -- (2,0) -- (2,1);
\end{tikzpicture} \\
Let vertex $v \in A$ be adjacent to vertex $s_2 \in T$.
We claim that the second neighbor of $v$ in $T$ can not be $s_4,s_5,s_6$ or $s_7$.
Assume by contradiction that vertex $v$ is adjacent to vertex $s_4$. We remove vertex $s_4$ from $S$ and add vertex $v$ to $S$, thus getting
a contradiction to condition $(5)$. By symmetry the same argument holds for vertices $s_5,s_6,s_7$.
Now as vertex $v$ is not adjacent to $s_4,s_5,s_6$ or $s_7$.
we get a contradiction to Claim $3$ for vertex $s_2$. \\
\noindent\textbf{Case 7.7: Tree $T$ is the following tree}. \\
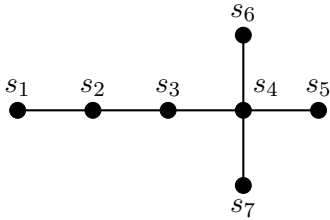
\begin{tikzpicture}
\draw[fill=black] (0,0) circle (3pt);
\draw[fill=black] (1,0) circle (3pt);
\draw[fill=black] (2,0) circle (3pt);
\draw[fill=black] (3,0) circle (3pt);
\draw[fill=black] (4,0) circle (3pt);
\draw[fill=black] (3,1) circle (3pt);
\draw[fill=black] (3,-1) circle (3pt);
\node at (0,0.3) {$s_1$};
\node at (1,0.3) {$s_2$};
\node at (2,0.3) {$s_3$};
\node at (3.3,0.3) {$s_4$};
\node at (4,0.3) {$s_5$};
\node at (3,1.3) {$s_6$};
\node at (3,-1.3) {$s_7$};

\draw[thick] (0,0) -- (1,0) -- (2,0) -- (3,0) -- (4,0);
\draw[thick] (3,-1) -- (3,0) -- (3,1);
\end{tikzpicture} \\
Let vertex $v \in A$ be adjacent to vertex $s_1 \in T$.
By Claim $4$ the second neighbor of $v$ in $T$ must be $s_5$,$s_6$ or $s_7$.
Assume w.l.o.g that vertex $v$ is adjacent to vertex $s_5$. We remove vertex $s_1$ from $S$ and add vertex $v$ to $S$, thus getting
a contradiction to condition $(5)$. \\
\noindent\textbf{Case 7.8: Tree $T$ is the following tree}. \\
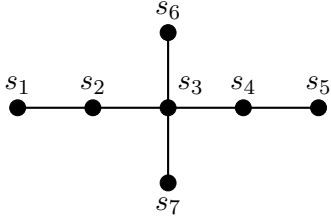
\begin{tikzpicture}
\draw[fill=black] (0,0) circle (3pt);
\draw[fill=black] (1,0) circle (3pt);
\draw[fill=black] (2,0) circle (3pt);
\draw[fill=black] (3,0) circle (3pt);
\draw[fill=black] (4,0) circle (3pt);
\draw[fill=black] (2,1) circle (3pt);
\draw[fill=black] (2,-1) circle (3pt);
\node at (0,0.3) {$s_1$};
\node at (1,0.3) {$s_2$};
\node at (2.3,0.3) {$s_3$};
\node at (3,0.3) {$s_4$};
\node at (4,0.3) {$s_5$};
\node at (2,1.3) {$s_6$};
\node at (2,-1.3) {$s_7$};

\draw[thick] (0,0) -- (1,0) -- (2,0) -- (3,0) -- (4,0);
\draw[thick] (2,-1) -- (2,0) -- (2,1);
\end{tikzpicture} \\
Let vertex $v \in A$ be adjacent to vertex $s_3 \in T$.
By Claim $4$ the second neighbor of $v$ in $T$ can not be $s_1$ or $s_5$.
Hence we get a contradiction to Claim $3$ for vertex $s_3$. \\
\noindent\textbf{Case 7.9: Tree $T$ is the following tree}. \\
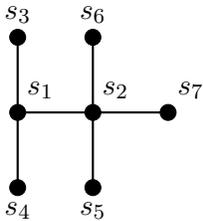
\begin{tikzpicture}
\draw[fill=black] (2,0) circle (3pt);
\draw[fill=black] (0,0) circle (3pt);
\draw[fill=black] (1,0) circle (3pt);
\draw[fill=black] (0,1) circle (3pt);
\draw[fill=black] (0,-1) circle (3pt);
\draw[fill=black] (1,-1) circle (3pt);
\draw[fill=black] (1,1) circle (3pt);
\node at (2.3,0.3) {$s_7$};
\node at (0.3,0.3) {$s_1$};
\node at (1.3,0.3) {$s_2$};
\node at (0,1.3) {$s_3$};
\node at (0,-1.3) {$s_4$};
\node at (1,-1.3) {$s_5$};
\node at (1,1.3) {$s_6$};

\draw[thick] (0,0) -- (1,0) -- (2,0);
\draw[thick] (0,-1) -- (0,0) -- (0,1);
\draw[thick] (1,-1) -- (1,0) -- (1,1);
\end{tikzpicture} \\
Let vertex $v \in A$ be adjacent to vertex $s_2 \in T$.
We claim that the second neighbor of $v$ can not be $s_3$ or $s_4$.
Assume that $v$ is adjacent to $s_3$. We can remove $s_3$ from $S$ and add $v$ to $S$, thus getting a contradiction
to condition $(4)$. In the same manner we can show that the second neighbor of $v$ in $T$ can not be $s_4$.
Hence we get a contradiction to Claim $3$ for vertex $s_2$. \\
\noindent\textbf{Case 7.10: Tree $T$ is the following tree}. \\
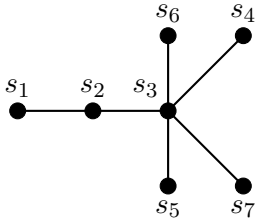
\begin{tikzpicture}
\draw[fill=black] (0,0) circle (3pt);
\draw[fill=black] (1,0) circle (3pt);
\draw[fill=black] (2,0) circle (3pt);
\draw[fill=black] (3,-1) circle (3pt);
\draw[fill=black] (3, 1) circle (3pt);
\draw[fill=black] (2,-1) circle (3pt);
\draw[fill=black] (2,1) circle (3pt);
\node at (0,0.3) {$s_1$};
\node at (1,0.3) {$s_2$};
\node at (1.7,0.3) {$s_3$};
\node at (3,1.3) {$s_4$};
\node at (3,-1.3) {$s_7$};
\node at (2,-1.3) {$s_5$};
\node at (2,1.3) {$s_6$};

\draw[thick] (0,0) -- (1,0) -- (2,0);
\draw[thick] (2,0) -- (3,1);
\draw[thick] (2,0) -- (3,-1);
\draw[thick] (2,0) -- (2,1);
\draw[thick] (2,0) -- (2,-1);
\end{tikzpicture} \\
Let vertex $v \in A$ be adjacent to vertex $s_3 \in T$.
By Claim $4$ the second neighbor of $v$ can not be $s_1$.
Hence we get a contradiction to Claim $3$ for vertex $s_3$.

\section{Appendix C}\label{last_appendix3}
It was shown in \cite{DBLP:journals/jgt/AlonMT01} that for any graph $G$ of maximum degree $4$ we have $a(G) \geq \frac{n}{2}$.
In this section we will improve this bound slightly by proving the following theorem.
\BTHM\label{triangle_free_deg_4}
Let $G=(V,E)$ be a triangle-free graph of order $n$ and average degree at most $4$. Then $a(G) \geq \frac{15n}{29}$.
\ETHM
We note that Example $2.2$ in \cite{DBLP:journals/jgt/AlonMT01} shows a $4$-regular graph $G$ on $n=14$ vertices for which $a(G) = \frac{4n}{7}$. \\
First we shall prove the following lemma.
\BL\label{simple_small_bound11}
Let $G=(V,E)$ be a triangle-free graph of order $n$ and average degree at most $4$. Then $a(G) \geq \frac{n+1}{2}$
\EL
\BPF
Assume w.l.o.g that graph $G$ is connected. By Theorem \ref{fundumental_bound_triangle_free_chinese} we have
\begin{equation}\label{good_bound_chinese1}
a(G) \geq \frac{10n - 5}{19}.
\end{equation}
Hence for $n>10$ we have $a(G) \geq \frac{10n - 5}{19} > \frac{n}{2}$.
The remaining case is when $n \leq 10$. If graph $G$ contains a vertex of degree at least $5$ then we are done as this vertex and $5$ of its neighbors are a tree of size $6$
 (as $G$ is triangle-free).
Hence we may assume that $G$ is of maximum degree $4$. 
Now if graph $G$ is not $4$-regular then $a(G) > \frac{n}{2}$ by Theorem \ref{fundumental_bound_triangle_free_alon}.
Thus we may assume that $G$ is a $4$-regular graph on at most $10$ vertices. \\
If $n \leq 9$ then we pick an arbitrary vertex $v$ in $G$ and its four neighbors thus getting a tree on $5$ vertices (as $G$ is triangle-free) and we are done. \\
The remaining case is when $G$ is a triangle-freen, $4$-regular graph on exactly $10$ vertices. We will assume that $a(G) \leq 5$ and get a contradiction. \\
Let $v_1 \in G$ be an arbitrary vertex and let $A=\{v_2,v_3,v_4,v_5\}$ be the set of neighbors of $v_1$ in $G$. Let $B=A \cup \{v_1\}$. Notice that as graph $G$ is triangle-free, set $B$ induces a tree in $G$. Let $C = V \backslash B$.  Since $|B|=5$ each vertex in $C$ must have at least $2$ neighbors in $A$ (otherwise we will get a forest on $6$ vertices in $G$).
Hence there must be at least $3$ vertices $v_6,v_7,v_8$ in $C$ each with exactly $2$ neighbors in $A$ for there are $12$ edges between the sets $A$ and $C$.
This means that at least two vertices from vertices $v_6,v_7,v_8$ are adjacent to the same vertex in $A$. Assume w.l.o.g that vertices $v_6,v_7$ are adjacent to vertex $v_2$.
As graph $G$ is triangle-free vertices $v_6$ and $v_7$ are not adjacent. Hence the set $(B \cup \{v_6,v_7\}) \backslash \{ v_2 \}$ induces a forest of size $6$ in $G$.
We got a contradiction and thus we are done.
\EPF \text{} \\
\noindent\textbf{Proof of Theorem \ref{triangle_free_deg_4}}: \\
Let $G'$ be a connected component of $G$ on $n'$ vertices.
By Lemma \ref{simple_small_bound11} and Equation \ref{good_bound_chinese1} we have
\begin{equation}\label{double_bound17}
a(G') \geq \max \left( \frac{n'+1}{2} , \frac{10n' - 5}{19} \right)
\end{equation}
Now notice the following.
\begin{itemize}
  \item For $n' \leq 29$ we have $\frac{15n'}{29} \leq \frac{n'+1}{2}$
  \item For $n' \geq 29$ we have  $\frac{15n'}{29} \leq \frac{10n' - 5}{19}$
\end{itemize}
We conclude by the observation above and Inequality \ref{double_bound17} that
\[
a(G') \geq \frac{15n'}{29}
\]
and as this holds for any connected component $G'$ of $G$ the theorem follows.
\qed

\section*{Acknowledgements}
Work supported in part by the Israel Science Foundation (grant No. 1388/16).

\bibliographystyle{alpha}




\end{document}